\documentclass{article}
\usepackage{amsmath}
\usepackage{booktabs}
\usepackage{dsfont}
\usepackage{xcolor}
\usepackage{graphicx} 
\usepackage{subcaption}
\usepackage{geometry}
\usepackage{authblk}
\usepackage{multirow}
\usepackage[capitalise]{cleveref}
\usepackage{float}
\DeclareMathOperator*{\argmin}{arg\,min}

\title{Solving QUBOs with a quantum-amenable branch and bound method}
\author[1]{Thomas H\"aner\thanks{thaener@amazon.com}}
\author[1]{Kyle E. C. Booth}
\author[2]{Sima~E.~Borujeni}
\author[2]{Elton~Yechao~Zhu}
\affil[1]{Amazon Quantum Solutions Lab, Seattle, WA 98170, USA}
\affil[2]{Fidelity Center for Applied Technology, FMR LLC, Boston, MA 02210, USA}
\date{}

\begin{document}

\maketitle

\begin{abstract}
Due to the expected disparity in quantum vs. classical clock speeds, quantum advantage for branch and bound algorithms is more likely achievable in settings involving large search trees and low operator evaluation costs. Therefore, in this paper, we describe and experimentally validate an exact classical branch and bound solver for quadratic unconstrained binary optimization (QUBO) problems that matches these criteria. Our solver leverages cheap-to-implement bounds from the literature previously proposed for Ising models, including that of Hartwig, Daske, and Kobe from 1984. We detail a variety of techniques from high-performance computing and operations research used to boost solver performance, including a global variable reordering heuristic, a primal heuristic based on simulated annealing, and a truncated computation of the recursive bound. We also outline a number of simple and inexpensive bound extrapolation techniques. Finally, we conduct an extensive empirical analysis of our solver, comparing its performance to state-of-the-art QUBO and MaxCut solvers, and discuss the challenges of a speedup via quantum branch and bound beyond those faced by any quadratic quantum speedup.\let\thefootnote\relax\footnotetext{Fidelity Public Information}
\end{abstract}

\section{Introduction}

Quantum algorithms have been an active area of research in chemistry \cite{lee2021even,von2021quantum}, machine learning \cite{schuld2015introduction}, and optimization \cite{farhi2014quantum,montanaro2020quantum}. Within optimization, there have been many quantum algorithms proposed, including those for more efficiently solving linear programs \cite{casares2020quantum,nannicini2024fast}, approximately solving combinatorial optimization problems \cite{farhi2014quantum}, and those for exactly solving combinatorial optimization problems \cite{montanaro2020quantum,booth2021quantum,chakrabarti2022universal}. Algorithms for the latter usually concern accelerating tree search schemes, often, specifically, branch and bound.

Branch and bound (B\&B), initially proposed for linear programming with discrete variables \cite{land2010automatic}, is an approach for solving optimization problems in a divide-and-conquer fashion. The algorithm (usually) avoids a worst-case complete enumeration search through the application of a bounding mechanism on the objective function that is used for deducing that a given subproblem need not be further explored as it, provably, cannot contain the optimal solution. B\&B is used ubiquitously across many industries and provides the foundation for many state-of-the-art commercial solvers, including Gurobi \cite{gurobi2021gurobi}, IBM CPLEX \cite{cplex2009v12}, and FICO Xpress \cite{berthold2018parallelization}. These solvers typically combine B\&B with cutting planes in a routine termed branch and cut to exactly solve large problems containing many variables and constraints \cite{bixby2012brief}. 

In recent years, there has been a specific focus in the quantum literature on quantum algorithms for B\&B \cite{montanaro2020quantum, chakrabarti2022universal}. In this setting, B\&B is modeled as an algorithm with access to `cost' and `branch' operators which, together, define a search tree. The cost operator provides provable lower (resp. upper) bounds on the objective function for minimization (resp. maximization) problems, and the branch operator divides the current (sub)problem into further subproblems. The proposed quantum branch and bound (QBB) algorithms tout a near-quadratic reduction, over their classical counterparts, in the number of times these operators are called (i.e., query complexity). The algorithms themselves leverage existing quantum-algorithmic routines for tree search (e.g., \cite{Belovs2013op}) and tree size estimation (e.g., \cite{Ambainis2017}). 

The analysis of QBB performance in prior work has focused on asymptotic query complexity and does not consider the magnitude of the constants that appear in the runtime of the oracle operators, which can be quite large. For the claimed speedups to hold in practice, assuming the oracle operations are the same between the classical and quantum methods, the quadratic reduction in tree size (i.e., calls to the cost and branch operators) has to neutralize the slower clock speed of the (fault-tolerant) quantum computer. Recent estimates \cite{babbush2021focus} find a relative quantum fault-tolerant vs. classical clock speed of $6.3\times 10^7$, even with the favorable assumption that a classical computer performs a single AND gate per cycle---an extremely optimistic scenario given the availability of classical multi- and many-core systems with low-latency Single Instruction Multiple Data (SIMD) floating-point and integer instructions. In this quantum-optimistic case, ignoring all further overheads that occur in an implementation of QBB\footnote{Overheads include, among others, (1) factors of $\sqrt d$ or $d^{1.5}$, (2) log-factors due to the binary searches in QBB, (3) overheads due to a reversible implementation of the required oracles.}, the tree size would have to be greater than $(6.3\times 10^7)^2 \approx 4.0\times 10^{15}$ \cite{babbush2021focus} before speedups could be realized. 

Based on the above, we can conclude that quantum advantage for B\&B is more likely when large trees are involved. Large trees, however, require a large number of calls to the cost and branch operators; for this reason, the cost of evaluating these operators should be made as low as possible. A common cost operator used in B\&B for integer programming involves solving a linear programming relaxation of the original problem, using methods such as simplex or interior point. This is prohibitively expensive to run on a quantum computer. Recent estimates of the expected cost of a quantum interior point method~\cite{PRXQuantum.4.040325} suggest that running such algorithms on a quantum computer at application scale is infeasible, at least with current algorithmic and quantum error correction techniques. Consequently, running a quantum interior point method as a subroutine of QBB is also impractical. 

Other optimization problems, however, may admit more practical cost operators. One such problem is the binary unconstrained quadratic optimization (QUBO) problem, where the binary nature of all variables and the absence of constraints yield a setting that may admit cheap cost operators. Previous work by Montanaro~\cite{montanaro2020quantum} identified Ising spin systems as a potential candidate for QBB, and it is well-known that Ising models and QUBO models are equivalent (see Section 2). Recognizing that the success of QBB is more likely in high search (i.e., large tree) settings and for problems that admit cheap effective bounding operators, we explore the development of a classical B\&B-based solver for QUBO problems.

Our contributions are as follows:
\begin{itemize}
    \item We design and implement a classical B\&B-based solver for QUBO problems, leveraging previously existing bounds proposed for Ising models. In this paper, we detail the main design considerations of this solver. 
    \item We conduct extensive empirical analysis on a variety of benchmarks (QUBO and Ising/MaxCut) and demonstrate the efficacy of our solver versus commercial integer programming solvers. Specifically, we show that our solver outperforms these commercial solvers for certain instance classes, setting the stage for potential quantization (and speedup) in the future.
    \item We conclude, however, that more research and optimizations are needed for low-inference solvers to become competitive with the state of the art. As long as low-inference solvers are outperformed by their high-inference counterparts, QBB is unlikely to yield a quantum advantage in practice, even with significant improvements in quantum computing hardware and software components.
\end{itemize}

The remainder of this paper is structured as follows. In Section \ref{sec:background} we provide necessary background. In Section \ref{sec:implementation} we describe the specific implementation of our solver, and in Section \ref{sec:experimental} we present the results of our experimental analysis. Finally, in Section \ref{sec:conclusions} we provide concluding remarks and future research directions. 

\section{Background} \label{sec:background}

In this section, we provide details on the fundamental concepts used within our cheap-inference B\&B-based solver for QUBO problems. 

\subsection{QUBO, Ising, and MaxCut models}
A QUBO model can be used to express an optimization problem that takes the following form: 
\begin{equation}
\min \sum_{i,j} Q_{i,j}x_ix_j 
\end{equation}
where $x_i \in \{0,1\}$ is a binary variable and $Q$ is a square matrix. 

A closely related problem is that of spin glasses, where an Ising model to minimize system energy is expressed as follows:
\begin{equation}\label{eq:ising-model}
   \min \sum_{i,j} J_{i,j} s_i s_j
\end{equation}
where $s_i \in \{-1, +1\}$ represents `up' and `down' spins, and $J$ is a matrix of interactions between pairs of spins. A linear term $\sum_i h_i s_i$ is often added to the Ising model to represent the presence of external fields. 

The Ising model can be converted into a QUBO through a simple substitution, e.g., $s_i \leftarrow 1 - 2x_i$. As a result, approaches developed for solving QUBOs can be readily applied to Ising problems. Conversely, approaches for solving Ising problems with external fields can be applied to QUBO problems using the reverse substitution, i.e., $x_i \leftarrow \frac{1-s_i}2$.

Minimizing the energy of an Ising model without external fields (as in \cref{eq:ising-model}) is equivalent to finding a maximum cut (MaxCut) of the spin coupling graph, i.e., the graph where nodes $i$ and $j$ are connected by an edge with weight $J_{i,j}$.

In the remainder of this section, we introduce two different bounds from the literature that have been previously applied in B\&B approaches to Ising model optimization \cite{kobe1978exact,hartwig1984recursive}. Both bounds can be implemented using relatively few elementary gates, and B\&B with these operators (usually) features large trees.

\subsection{Kobe-Hartwig bound}

In 1978, Kobe and Hartwig proposed a B\&B method for computing an exact ground state of Ising systems~\cite{kobe1978exact}, i.e., for computing $s^\star\in\{-1,1\}^n$ such that
\begin{equation}
    s^\star \in \argmin_{s\in\{-1,1\}^n} \sum_{i<j} J_{i,j} s_i s_j,
\end{equation}
where $J_{i,j}$ is the real-valued coupling strength between spins $s_i$ and $s_j$.

Their method starts with a root node bound $B_\text{KH}(\star)$ given by
\begin{equation}\label{eq:kh-root}
    B_\text{KH}(\star)=-\sum_{i<j}|J_{i,j}|.
\end{equation}
As individual spins are assigned a value during branching, the bound is strengthened with a simple update. Assuming that $s_1$ up to $s_k$ have already been assigned, the updated bound after branching on $s_{k+1}$ can be computed via
\begin{equation}\label{eq:kh}
    B_\text{KH}(s_1, ..., s_{k+1}) = B_\text{KH}(s_1, ..., s_k) + 2\sum_{i\leq k} |J_{i,k+1}|\left[J_{i,k+1}s_is_{k+1} > 0\right],
\end{equation}
where $[\cdot]$ denotes the Iverson bracket and $B_\text{KH}(s_1,...,s_k)$ denotes the Kobe-Hartwig bound for the partial assignment where $s_1,...,s_k$ have been assigned and $s_{k+1},...,s_n$ have not.
\subsection{Hartwig-Daske-Kobe bound}

In follow-up work to Ref.~\cite{kobe1978exact}, Hartwig, Daske, and Kobe propose an improved bound~\cite{hartwig1984recursive}. In contrast to the Kobe-Hartwig bound, this new bound is based on the optimal solution, $E_l$, of lower-dimensional subproblems with $l=2, ..., (n-1)$ spin variables.

At a B\&B tree node where spin variables $s_1,...,s_k$ have already been assigned, the Hartwig-Daske-Kobe bound $B_\text{HDK}$ is computed as follows
\begin{equation}\label{eq:hdk}
    B_\text{HDK}(s_1,...,s_k)=\sum_{i<j\leq k}J_{i,j}s_is_j - \sum_{j>k}\left|\sum_{i\leq k} J_{i,j}s_i\right| + E_{n-k},
\end{equation}
where $E_{n-k}$ is the optimal solution to the $(n-k)$-dimensional subproblem that is obtained by canceling the first $k$ rows and columns of the coupling matrix $J$. When branching on the first spin variable $s_1$, one needs access to the $(n-1)$-dimensional subproblem solution $E_{n-1}$. In turn, solving this subproblem via B\&B requires $E_{n-2}$ to branch on the first variable. Therefore, one first solves all lower-dimensional problems starting with a two-dimensional problem, the solution to which is $E_2$, all the way up to an $(n-1)$-dimensional problem, from which one obtains $E_{n-1}$.

If the order in which variables are assigned is fixed, then these lower-dimensional subproblems have to be solved only once. While there is a tradeoff between the potential benefits of variable reordering and the required recomputation of lower-dimensional subproblems, we opt for a fixed variable ordering in our implementation.

As was the case for the Kobe-Hartwig bound, it is possible to obtain an analogous bound that is valid for objective functions with an additional field term. We also note that, instead of the optimal solutions $E_k$, any lower bound to the true value $E_k$ would also result in a correct B\&B solver. In particular, we will use this fact to combine $B_\text{HDK}$ with $B_\text{KH}$ in order to reduce the time spent solving lower-dimensional subproblems.
\subsection{Inference costs and tree sizes}

In this paper, we focus on B\&B methods with cheap inference, as measured, for example, by the number of binary AND and exclusive-OR gates required to implement the inference step at a given problem size. At the same time, the corresponding B\&B tree (and, thus, the quantum speedup from QBB) should be large, but not too large, such that the solver still performs well.

Computing the Kobe-Hartwig bound is possible using just $\mathcal O(n)$ additions for each parent-to-child transition, \cref{eq:kh}, resulting in a total cost per leaf of $\mathcal O(n^2)$ additions (assuming caching of parent bounds). While this cost is exceptionally low, it is also not a particularly effective bound. As a result, solving Ising and QUBO problems with significantly more than $40$ variables is impractical using the KH bound.

\begin{figure}[ht]
    \centering
    \begin{subfigure}{0.49\linewidth}
    \includegraphics[width=\linewidth]{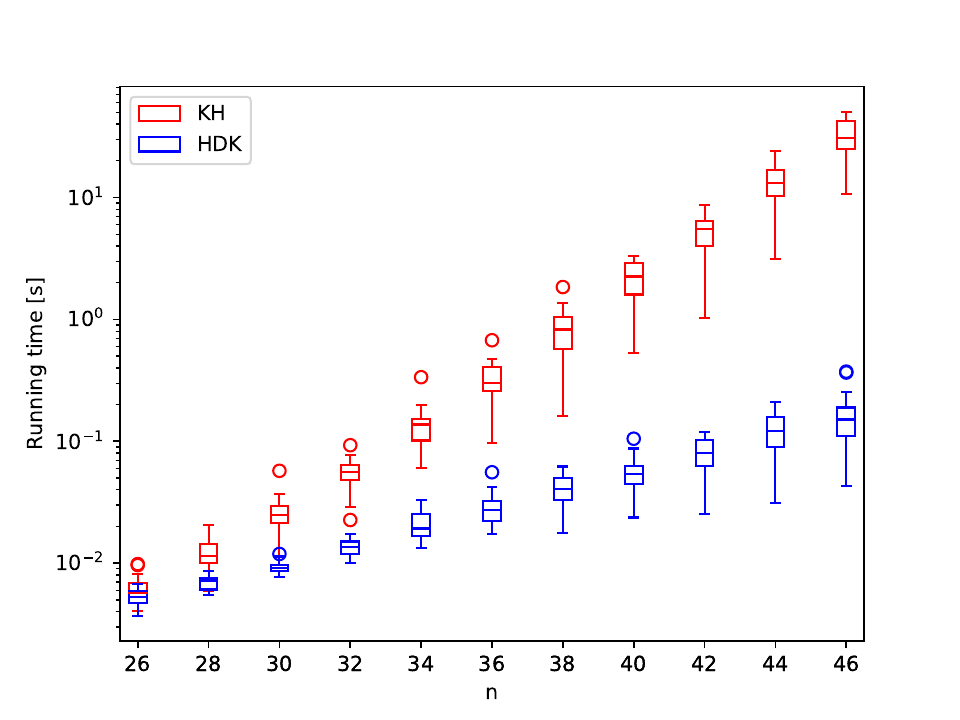}
    \caption{Runtime comparison}
    \label{fig:kh_vs_hdk_time}
\end{subfigure}
\begin{subfigure}{0.49\linewidth}
    \includegraphics[width=\linewidth]{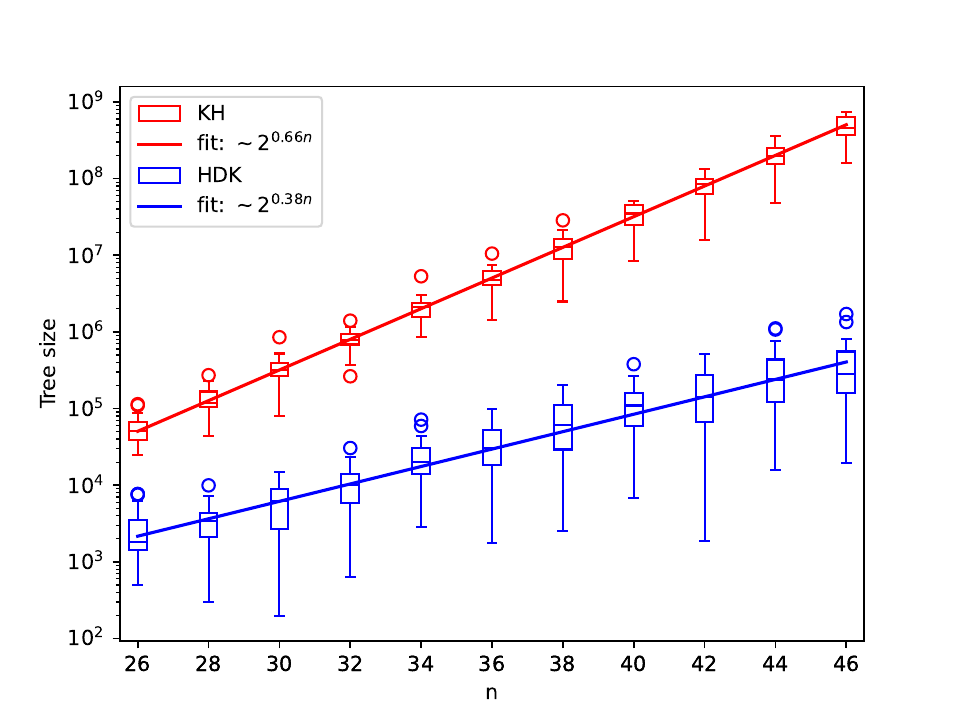}
    \caption{Tree size comparison}
    \label{fig:kh_vs_hdk_size}
\end{subfigure}
\caption{Comparison between an improved version of the Kobe-Hartwig (KH) bound~\cite{kobe1978exact} (see text) and the Hartwig-Daske-Kobe (HDK) bound~\cite{hartwig1984recursive}. Each box represents the single-threaded runtime (in \ref{fig:kh_vs_hdk_time}) or tree size (in \ref{fig:kh_vs_hdk_size}) across 20 Ising instances with $n$ spin variables, where the integer couplings $J_{i,j}$ and fields $h_i$ were chosen uniformly at random from $(-100,100)$. While, compared to the HDK bound, the KH bound is cheaper to evaluate and produces larger B\&B trees, which is advantageous for quantum branch-and-bound, it is too weak to yield a competitive solver.}
\label{fig:kh_vs_hdk}
\end{figure}

Even an improved version of the KH bound, where we borrow the second term from the HDK bound, \cref{eq:hdk}, and use the initial KH bounds for $E_{n-k}$, is only capable of solving Ising and QUBO instances with up to about 50 variables, see Fig.~\ref{fig:kh_vs_hdk}.

In contrast, the HDK bound is a much stronger bound, but it requires the solutions of $(n-k)$-dimensional subproblems for all $k\in\{1,...,n-2\}$. At first glance, this may seem excessively expensive. However, assuming that i) the solver runtime grows exponentially with $n$, ii) we are solving the problem to proven optimality, and iii) the ordering of branching variables is fixed, these lower-dimensional problems only contribute a (possibly large) constant factor to the total runtime. Consequently, under these assumptions, the HDK bound has a per-leaf cost of $\mathcal O(n^3)$ additions. Assuming perfect caching of parent bounds (e.g., depth-first search) and using $\mathcal O(n\cdot b)$ memory, where $b$ denotes the number of bits used to represent bound values and energies, parent-to-child updates can be achieved using only $\mathcal O(n)$ additions. As a result, the cost per leaf is $\mathcal O(n^2)$ $b$-bit additions in that case.

In terms of their asymptotic costs, the two bounds are surprisingly similar. Moreover, the constant factor due to lower-dimensional subproblems in the HDK bound may be reduced significantly by only computing subproblems with $k > k_\text{min}$ and using the KH bound for $k\in\{1,...,k_\text{min}\}$.

Finally, previous work~\cite{montanaro2020quantum} has found empirically that the HDK bound produces large trees with $\sim 2^{0.37n}$ nodes for $n$-spin Sherrington-Kirkpatrick instances and that it is possible to solve 50-spin instances within minutes on a standard laptop using this bound. While our results agree in terms of the B\&B tree sizes (Fig.~\ref{fig:kh_vs_hdk}) we find that our implementation of the HDK bound, together with various optimizations such as variable reordering, is capable of solving similar 50-spin instances in less than one second. Additionally, some instances from the BiqMac Library~\cite{wiegele2007biq} with more than 100 variables can be solved in less than an hour, as detailed in Section~\ref{sec:experimental}.

Therefore, the HDK bound---with its low cost, reasonable bound strength, and relatively large B\&B trees---is a promising candidate for a speedup via QBB. Next, we investigate whether the HDK bound is strong enough to make it a worthwhile candidate for a future speedup via QBB, i.e., whether it yields a classical B\&B solver that is sufficiently performant such that even small improvements in runtime through QBB would suffice for it to start outperforming state-of-the-art classical solvers.

\section{Implementation} \label{sec:implementation}

In this section, we describe the implementation of our solver, including primal/dual bounding mechanisms, tree traversal logic, variable and value ordering heuristics, and our multithreading approach.

\subsection{Bounding mechanisms}

\subsubsection{Primal bounds}\label{sec:primalheur}

To obtain primal bounds (i.e., feasible solutions that are not necessarily optimal), we use a combination of two methods. The first is a greedy extension of a lower-dimensional solution to a solution of a larger problem: Starting with a lower-dimensional solution, we pick the unassigned spin that, when assigned to $\{1,-1\}$, results in the largest change in energy, and we assign it the value that minimizes the energy. We repeat this process until all unassigned spins have been assigned.

The second primal heuristic is simulated annealing~\cite{kirkpatrick1983optimization}, where we initialize the search with either a random configuration or the configuration obtained from the greedy extension above. Simulated annealing is run once on the full problem before the solver starts solving the first subproblem. Moreover, for each subproblem, the solver uses simulated annealing to find a good solution to the subproblem before starting to explore the B\&B tree. Finally, simulated annealing (with a very short annealing schedule to minimize the overhead) is also used in our variable reordering heuristic as described in Section~\ref{sec:reordering}.

\subsubsection{Dual bounds}

\paragraph{Evaluation of the HDK bound.}

In contrast to the Kobe-Hartwig bound, the HDK bound does not permit a straightforward and efficient reuse of node bounds to compute the bounds of their children nodes. To enable such reuse, we store additional intermediate results at each node. Specifically, in addition to $\sum_{i<j,j\leq k} J_{i,j}s_is_j$, a node at level $k$ stores $\sigma_j=\sum_{i\leq k} J_{i,j}s_i$ for each $j>k$. Using these intermediate results, the second sum in \cref{eq:hdk} may be evaluated more efficiently: computing a child node bound from the parent bound is possible using $\mathcal O(n)$ operations and $\mathcal O(n)$ space instead of $\mathcal O(n^2)$ operations and $\mathcal O(1)$ space.

\paragraph{Combined HDK and KH bounding.}

The total time spent solving all lower-dimensional subproblems constitutes a significant fraction of the total runtime. Further, it was found that few nodes (if any) are truncated near root node. The bounds of these nodes involve $E_{n-k}$ for small $k$. These are the optimal solutions to the largest subproblems being solved in the precomputation phase, i.e., before the B\&B on the full $n$-dimensional problem is started.

We make use of these observations by only solving subproblems up to dimension $(n-k_\text{min}-1)$, where $k_\text{min}\sim 20$ was found empirically. We compute the bounds for the first $k_\text{min}$ levels by replacing $E_{n-k}$ in \cref{eq:hdk} by the Kobe-Hartwig root bound \cref{eq:kh-root}. 

\paragraph{Extrapolation from subproblems.}

For large and difficult problems, the precomputation phase can be long. In order to provide valid primal and dual bounds throughout the entire solver runtime, we simply extend both primal and dual bounds of subproblems to the full problem.

We arrive at valid primal (upper) bounds for the $n$-dimensional problem using our two primal heuristics (greedy extension and simulated annealing, see Section~\ref{sec:primalheur}). 

While solving for $E_{n-k}$, a valid lower bound $B_{n-k}$ on $E_{n-k}$ may be extrapolated to a valid lower bound $B_n^{k}$ for $E_n$ via
\begin{equation}
    E_n \geq B_n^k := -\sum_{i<j} |J_{i,j}| + \sum_{i\geq k,i<j} |J_{i,j}| + B_{n-k}.
\end{equation}

\paragraph{Local field term.} \label{sec:local-field}While an extension of the two bounds with an external field term $\sum_ih_is_i$ is straightforward, there are two options for the HDK bound. The first is to assign $h=0$ for all lower-dimensional subproblems and only include it in the final B\&B in the first two summands in \cref{eq:hdk}. The second option is to keep the local fields in lower-dimensional subproblems and only add it to the first term in \cref{eq:hdk}, which computes the energy of and among spins that have already been assigned.

The two options have different performance characteristics and future work could investigate possible heuristics for making an appropriate choice. In our benchmarking section, we simply opt for the second option, but we have found that some instances could be solved significantly faster using the first (and vice-versa).

There are several reasons for this behavior. Among others i) including local fields in subproblems can make them significantly cheaper to solve, for example, if fields are strong relative to the spin-couplings, and ii) adding local fields to the second term in \cref{eq:hdk} instead can lead to cancellations and, thus, a tighter bound that may allow the solver to eliminate more nodes.

\subsection{Tree traversal}

\paragraph{Depth-first traversal.} Given that all branching variables are $s_i\in\{-1,1\}$, we may visit all $n$-variable assignments by iteratively incrementing an integer $x$ from $x=0$ to $x=2^n-1$. For each $x$, the assignment of the $i$-th branching variable $s_i$ then corresponds to $(-1)^{x_i}$, where $x_i$ is the $i$-th bit in the binary representation of $x$. For example, the $3$-bit integer $x=6$ corresponds to the assignment $s=(1,-1,-1)$.

In B\&B, subtrees may be eliminated based on the bound, and an initially empty assignment is expanded, for example, one variable at a time. To this end, we store the depth of the current node, starting with $d=0$. The pair $(x,d)$ corresponds to a B\&B tree node with the assignment
\begin{equation}
    s_i\leftarrow a_i(x,d):=\left\{\begin{matrix} (-1)^{x_i},&\text{if $i<d$}\\\star,&\text{otherwise}\end{matrix}\right.,
\end{equation}
where $s_i=\star$ denotes that $s_i$ has not yet been assigned. Our depth-first traversal increments $d$ until it determines that the subtree may be eliminated (based on the bound). At this point, we add $2^{n-d}$ to $x$ (using inverted bit-significance for $x$) and update $d\leftarrow (n-k)$, where $k$ is the 0-based index of the first nonzero bit in the list $x_{n-1},x_{n-2},...,x_0$. The tree has been fully explored once $x$ is again equal to zero.

Given the low memory requirements and cheap update operation, the overhead of this B\&B tree data structure is very small, i.e., most of the solver running time is spent computing and updating the dual bound.

\paragraph{Dual bound and depth-first search.} While traversing the B\&B tree, our method keeps track of the best dual bound that is valid for all remaining (i.e., not yet visited) nodes. One possible choice is the bound of a tree node that dominates all remaining nodes. While efficiently computable and cheap to keep track of, the resulting bound remains loose for a significant portion of the B\&B tree exploration. In particular, the initial (root) bound may only be tightened once the entire subtree rooted at a child node of the root node has been traversed.

\paragraph{Dual bound and best-first search.}
In contrast to depth-first search (DFS), best-first search (BFS) keeps track of a node frontier. In each traversal step, it expands the node in the frontier with the highest potential to result in a good solution. Our method uses the node bound as an estimate for its potential. The resulting B\&B method uses more memory than the DFS approach, which, in turn, causes a slowdown of the solver. However, BFS also yields tighter dual bounds.

\paragraph{Hybrid BFS/DFS.} The final implementation of our solver uses a hybrid of the two approaches above, starting with a BFS until the frontier size exceeds a user-supplied limit. Whenever its size exceeds this limit, the entire subtree rooted at the next node in the frontier is traversed using DFS.

By increasing the maximal frontier size, one may tighten dual bounds at the expense of extra space and time overheads due to the larger node frontier.

\subsection{Global variable/value ordering}\label{sec:reordering}
While general B\&B solvers may choose the next variable to branch on using an arbitrary heuristic, our use of the HDK bound makes choosing different variables at a fixed distance from the root node much more costly. Choosing a different ordering of variables may invalidate all lower-dimensional bounds and the solver would have to recompute them recursively. As a result, we opt for the simplest solution and select a fixed global variable ordering instead.

Determining a fixed global variable ordering avoids the need for solving several instances of a given dimension $n_d<n$, and thus reduces the amount of time spent constructing the bound for the actual $n$-dimensional B\&B. However, we note that it could be beneficial in some cases to use non-fixed variable orderings, especially for nodes that are far from the root node, as long as the number of different orderings is modest. We leave investigating this tradeoff to future work.

Our variable reordering heuristic aims to combine two objectives. First, variables with strong couplings to other variables should end up toward the top of the branch-and-bound tree. Second, since we do not compute $E_{n-k}$ for small $k$, we may prioritize solving subproblems that we anticipate to be easy and trade these for $k_\text{min}$ 
difficult subproblems that end up at the very top of the final branch and bound tree. We thus allow the solver to reorder the last $2k_\text{min}$ variables according to a slightly different heuristic that attempts to combine the first with the second objective.

The two heuristics that we have found to work reasonably well are
\begin{align*}
    H_1(i)&=\sum_{j=1, j\neq i}^n |J_{i,j}| + |E_i(s^\star)|\\
    H_2(i)&=\sum_{j=1, j\neq i}^n |J_{i,j}| - E_i(s^\star),
\end{align*}
where $s^\star$ is a primal solution (found using simulated annealing) of the next subproblem and $E_i(s^\star)$ denotes the energy contribution from the $i$th row/column of the next subproblem.

When deciding which variable to place at position $n-i$ (corresponding to the subproblem of dimension $i$), the solver places a new variable at position $i$ and evaluates $H_1$ and $H_2$. This process is repeated for each of the variables that have not been added to the variable ordering. For positions $i$ from $1$ to $n-2k_\text{min}-1$, the solver picks the variable with the smallest $H_1$-value. For $i \geq n-2k_\text{min}$, the variable's $H_2$-value is used instead.

In addition to variable reordering, a B\&B solver may choose to branch on $s_i=1$ before considering the $s_i=-1$ case and vice versa. Aiming for a simple and efficient implementation, we also opt for a global value ordering based on the value of the corresponding external field $h_i$. If $h_i\geq 0$, our solver branches on $s_i=-1$ first, since then $h_is_i \leq h_i(s_i+2)$. Otherwise, it starts with $s_i=1$. 

\subsection{Multithreading}

To leverage multi- and many-core systems, our solver supports multithreading via OpenMP~\cite{dagum1998openmp}. Given a number of threads $n_t=2^{k_t}$ with $k_t\in\mathds N$, each thread with index $i_t\in\{0,...,n_t-1\}$ only visits nodes where the first $k_t=\log_2(n_t)$ bits of $x$ agree with those of $i_t$.

While a more elaborate multithreading scheme is likely to perform better (e.g., using OpenMP tasks), our approach is straightforward to implement and potential load balancing issues may be mitigated using hyperthreading.

To improve multithreading efficiency and thus allow the solver to better leverage systems with many cores, future work may implement alternative (e.g., task-based) strategies.

\section{Solver performance} \label{sec:experimental}

In this section, we present benchmarking results and compare the performance of our solver to that of Xpress, Gurobi, BiqMac, and BiqBin.

\subsection{Machine specifications}

We ran all benchmarks on two slightly different machines. Gurobi benchmarking was performed on a server equipped with an Intel Core Processor of the Haswell generation, 4 physical cores, 16GB of RAM and running at 2.9GHz, whereas Xpress and our solver were run on a 2.6GHz 6-core Intel Core i7 machine with 16GB of RAM running MacOS Sonoma 14.3. FICO Xpress v8.13.5 used up to 12 threads.

\subsection{Solver parameters}

\paragraph{Commercial solvers.} We implement QUBO and MaxCut models in the commercial solvers using the substitution described in Section \ref{sec:background}. Default solver parameters were used for both Xpress and Gurobi. The ``TimeLimit'' parameter was set to 3600 seconds for both solvers.

\paragraph{Our solver.} We set the maximum number of BFS nodes to $100\,000$ for all instances. For the bound precomputation, we chose to always keep local fields in subproblems. As discussed, it would be beneficial to remove local fields for some instances (e.g., gka$\star $b), and a heuristic would be needed to make this decision for each instance.

While our solver generally benefits from extensive hyperthreading (using 2-4 threads per core), since this remedies some load balancing issues that arise in our parallelization strategy, we only allowed up to $8$ threads on the $6$-core machine in order to facilitate a fair comparison with Xpress---which we ran using at most 12 threads. 

Finally, we enabled variable reordering for all instances except for the Ising set, for which a reordering heuristic would have to take into account the relative spin coupling strengths and connectivity.

\subsection{Benchmark instances}

For benchmarking, we use the instances from the BiqMac Library~\cite{wiegele2007biq}, which contains unconstrained binary quadratic programming, i.e. QUBO, as well as MaxCut instances, as introduced in Section~\ref{sec:background}.
For more information and additional details on each instance class, we refer the reader to Ref.~\cite{wiegele2007biq}.

The instances from the BiqMac Library have been used to benchmark various QUBO and MaxCut solvers, including BiqMac~\cite{rendl2010solving} and BiqBin~\cite{gusmeroli2022biqbin}. The BiqMac Library is thus a natural candidate for evaluating the performance of our solver.

\subsection{Comparison to Xpress and Gurobi}

\begin{figure*}[ht]
\includegraphics[width=\linewidth]{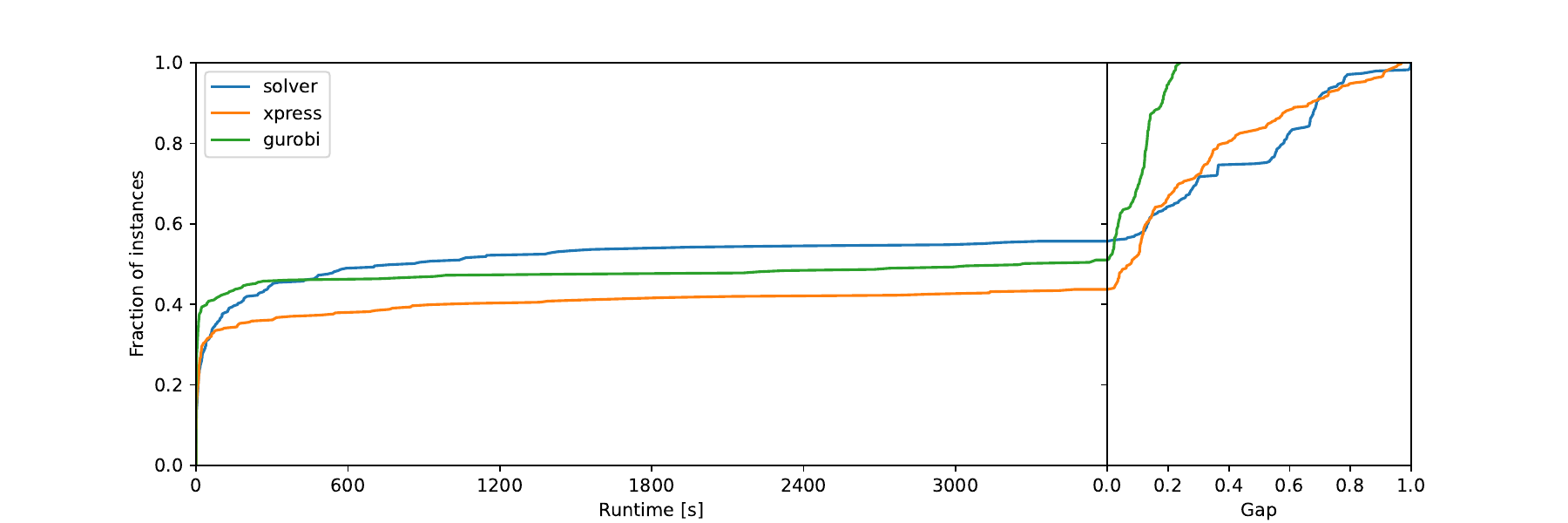}
\caption{Comparison summary of our solver, Xpress, and Gurobi for all 343 instances in biqmaclib. The plot depicts solver running time and the final gap if the time limit of 1h is reached.}
\label{fig:resall}
\end{figure*}

\begin{figure*}[ht]
\includegraphics[width=\linewidth]{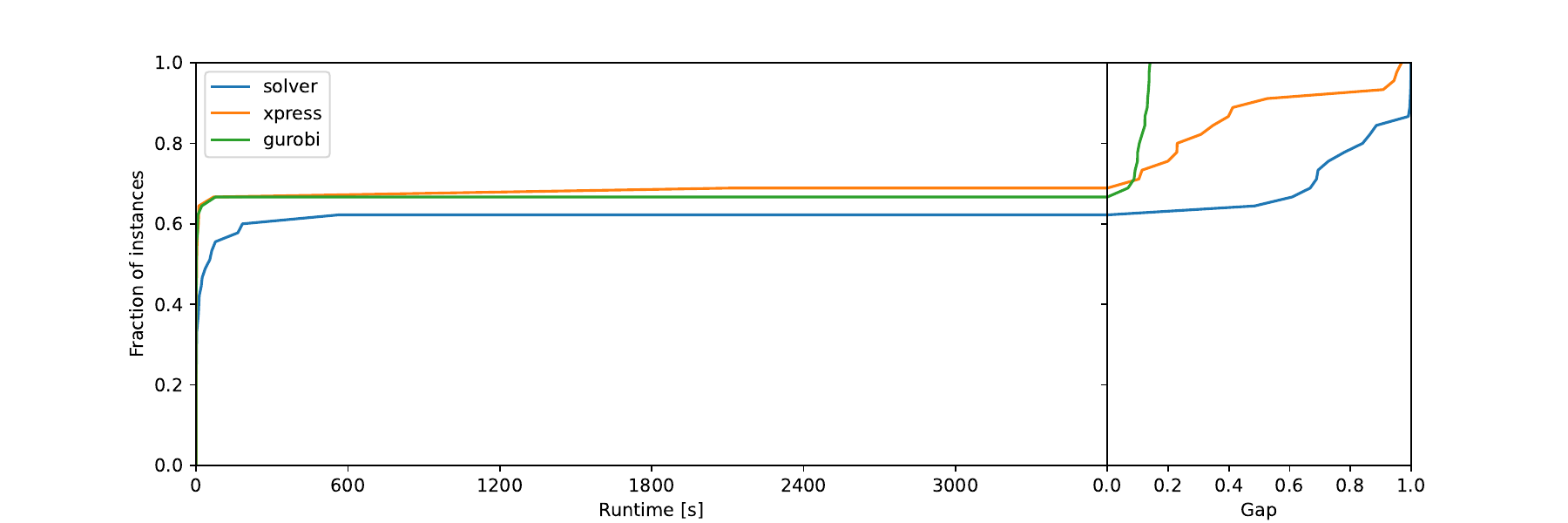}
\caption{Comparison of our solver, Xpress, and Gurobi for all 45 \texttt{gka*} instances in biqmaclib. The plot depicts solver running time and the final gap if the time limit of 1h is reached. Our solver performs worse than Gurobi and Xpress when keeping local fields in the bound precomputation phase. Removing local fields in lower-dimensional subproblems would allow our solver to solve more instances to optimality.}
\label{fig:resgka}
\end{figure*}

\begin{figure*}[ht]
\includegraphics[width=\linewidth]{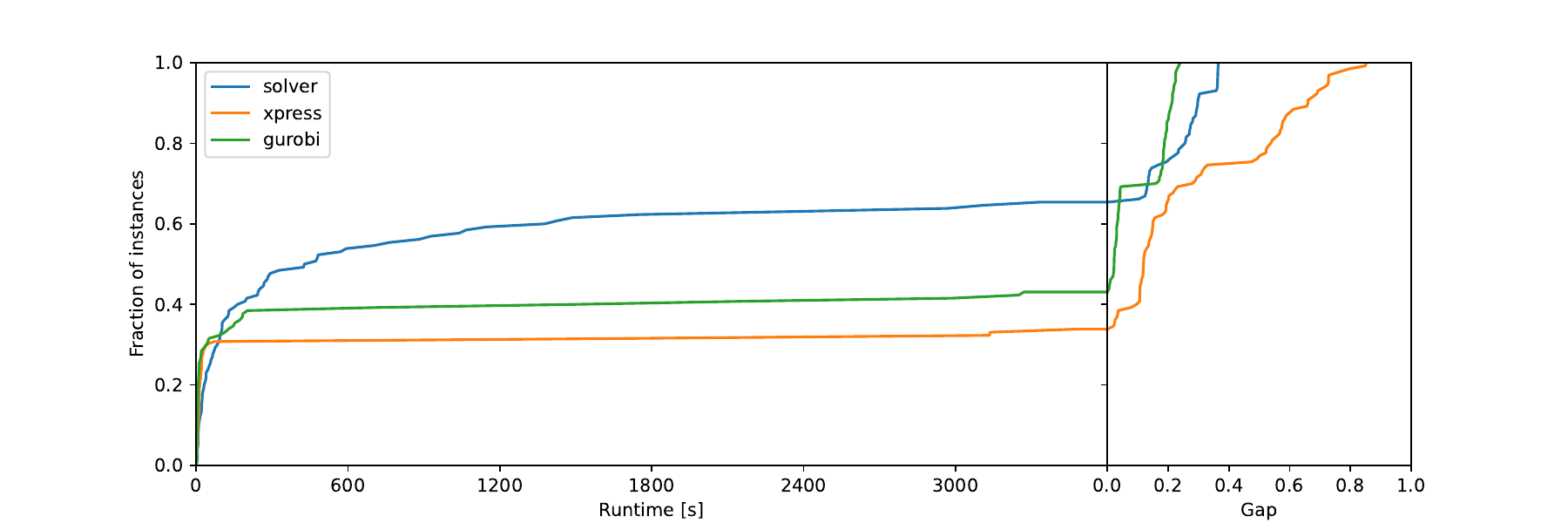}
\caption{Comparison of our solver, Xpress, and Gurobi for all 130 \texttt{rudy} instances in biqmaclib. The plot depicts solver running time and the final gap if the time limit of 1h is reached. The performance of our solver relative to LP-based solvers is best for dense problems, which allows it to outperform Gurobi and Xpress for the \texttt{rudy} test set.}
\label{fig:resrudy}
\end{figure*}

When comparing the performance of our solver to that of Gurobi and Xpress, we notice that our solver performs best for dense instances, such as \texttt{be100*}, \texttt{pm1d*}, and \texttt{w0*} from the BiqMac Library. Since there are many dense instances in the BiqMac Library, the overall comparison of solved instances looks favorable, see the left-hand side of Fig.~\ref{fig:resall} and Fig.~\ref{fig:resrudy}. However, we do find that Gurobi provides tighter bounds than our solver (and Xpress) for problems that are not solved to optimality within the time limit. This can be seen on the right-hand side of Fig.~\ref{fig:resall} and Fig.~\ref{fig:resgka}.

\begin{figure*}[ht]
\includegraphics[width=\linewidth]{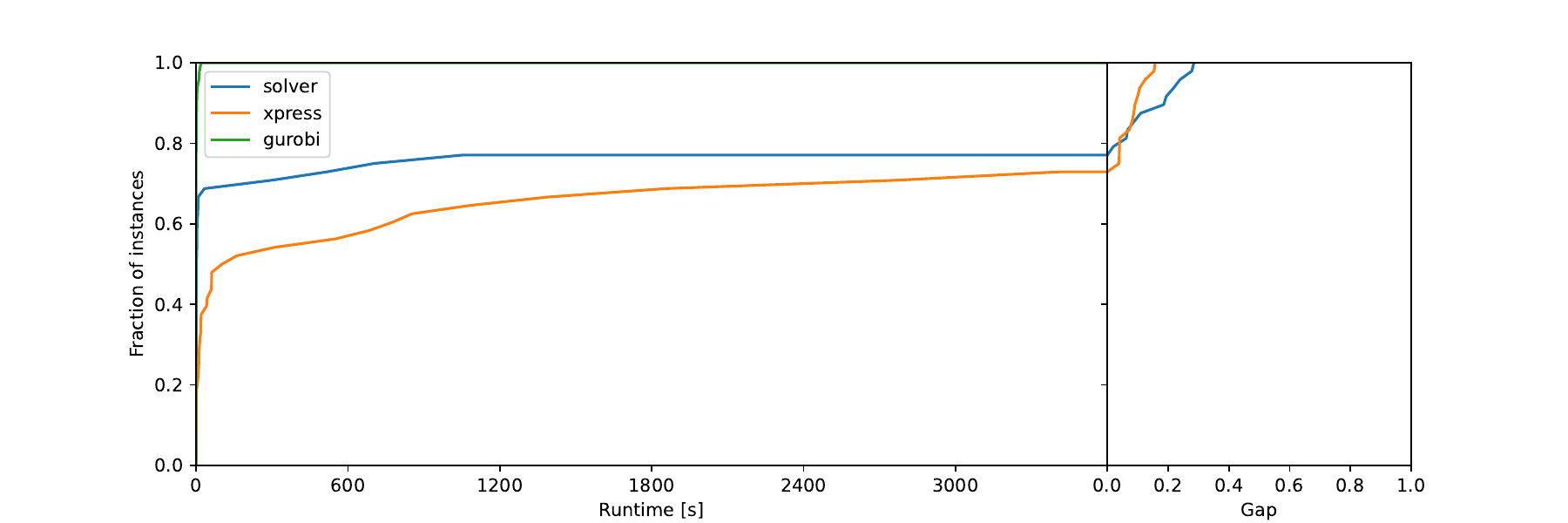}
\caption{Comparison summary of our solver, Xpress, and Gurobi for all 48 Ising instances in biqmaclib. The plot depicts solver running time and the final gap if the time limit of 1h is reached. Gurobi solves all instances in a few seconds.}
\label{fig:resising}
\end{figure*}

For sparse problems, LP-based methods clearly outperform our solver, as can be seen in Fig.~\ref{fig:resising}, which depicts a summary for various Ising instances, many of which are sparse (e.g., 2D and 3D grid connectivity). Gurobi solves all of these instances in a few seconds, whereas our solver fails to solve some of the larger 2D and 3D Ising instances within the time limit. In contrast, Xpress successfully solves all (sparse) 2D and 3D Ising instances, but fails to solve some of the dense instances, see Table~\ref{table:resultssummary_ising}.

\subsection{Comparison to BiqBin and BiqMac}

While Gurobi and Xpress perform well for sparse QUBO and MaxCut instances, there exist other solvers that outperform such LP-based methods by a wide margin, especially for dense instances. Two examples of such state-of-the-art QUBO and MaxCut solvers are BiqMac~\cite{rendl2010solving} and BiqBin~\cite{gusmeroli2022biqbin}.

These solvers have also been benchmarked using the instances from the BiqMac Library~\cite{wiegele2007biq}, allowing us to compare our solver to the state of the art for solving dense QUBO and MaxCut problems.
However, we note that a direct comparison to the numbers in Refs.~\cite{rendl2010solving, gusmeroli2022biqbin} is imperfect, especially in the case of BiqMac, which was evaluated over 10 years ago on a Pentium IV with 3.6 GHz and 2 GB of memory.

With this caveat in mind, we find that our solver performs quite well compared to BiqMac for some of the instance classes, e.g., $G_{-1/0/1}$, for which BiqMac requires on average 56 minutes and our solver solved eight out of 10 instances in less than 1 hour each. While there are other instance classes that are solved faster by our solver such as some of the Ising instances, we have seen that Gurobi is the solver of choice for these instances.

However, the newer and parallelized BiqBin solver~\cite{gusmeroli2022biqbin} clearly outperforms our solver: For example, even the sequential version of BiqBin is able to solve each of the \texttt{be100*} instances in less than 82 seconds. The maximum runtime of our solver for this set of instances is 545 seconds, significantly longer than BiqBin despite multithreading. Moreover, the advantage of BiqBin over our solver is even more pronounced for larger instances: While our solver only solves two \texttt{be120.8} instances in less than an hour, sequential BiqBin solves all of these instances with an average runtime of approximately two minutes (and max. runtime of 267 seconds).

We may thus conclude that additional optimizations and research would be required for this low-inference approach to be competitive with state-of-the-art QUBO and MaxCut solvers. As a consequence, also a future quantum advantage by applying QBB to the B\&B method considered in this work is unlikely without significant algorithmic improvements.

\section{Conclusion and future work} \label{sec:conclusions}

We found that it is possible to develop a B\&B based QUBO solver that performs reasonably well despite being on the high-exploration, low-inference end of the spectrum of B\&B solvers. However, we would like to see a performance that is (close to) on-par with or better than that of state-of-the-art QUBO solvers for QBB to yield a practical advantage on a future quantum computer. Thus, significant improvements are needed not only in quantum hardware and software (as is the case for any quadratic quantum speedup~\cite{babbush2021focus,hoefler2023disentangling}), but also in classical B\&B solvers to enable practical quantum speedups via QBB.

In addition to researching low-inference solvers for other application domains, future work may investigate various improvements to the implementation of our solver. A task-based multithreading/processing strategy may improve efficiency on multi-core systems and when distributing the computation across several compute nodes.
We list additional extensions below.

\paragraph{Heuristic for (non)zero fields in recursion.} Our solver supports two variants of the recursive bound: One where the field term is removed when solving subproblems, and one where field terms are kept (see Section~\ref{sec:local-field}). While our solver always used the latter in our benchmarking, using the former leads to substantial speedups in some cases (and slowdowns in others). For example, our solver timed out when solving \texttt{gka4b}, but successfully solves it to optimality in less than 30 seconds when field terms are removed in lower-dimensional subproblems. A heuristic may be used to decide which variant to use. One may even combine the two variants and thus allow the solver to make this decision on a per-variable basis.

\paragraph{Improved variable reordering heuristic.} The current variable ordering heuristic does not take into account spin-connectivity or, more generally, the relative coupling strengths among different sets of spins. However, this may be beneficial, e.g., for sparse problems and for problems where coupling strengths vary significantly.

In addition to an improved global variable ordering, future work could also investigate local reordering of variables: As long as the number of different orderings is small, the overhead from having to solve multiple subproblems of a given dimension may be compensated by a tighter dual bound.

\paragraph{Native support for simple constraints.} Our implementation does not support constraints natively, i.e., without introducing penalty terms. However, certain types of constraints could be supported by modifying the tree traversal logic. For example, certain cardinality constraints could be handled by only iterating over integer variables $x$, where the number of nonzero bits is equal to or less than a constant. This may improve the performance for subproblems (during the bound precomputation phase) as well as for the final B\&B.

\section{Acknowledgements}
This work is a collaboration between Fidelity Center for Applied Technology, Fidelity Labs, LLC., and the Amazon Quantum Solutions Lab. The Fidelity publishing approval number for this paper is 1155552.1.0.
\bibliographystyle{unsrt}

\section{Detailed results}

\subsection{Summary}

\begin{table}[H]\footnotesize
\centering
\begin{tabular}{ccccccccccc}
\multicolumn{2}{l}{\textbf{Parameters (\texttt{be*})}} & \multicolumn{3}{l}{\textbf{Solver}} & \multicolumn{3}{l}{\textbf{Xpress}} & \multicolumn{3}{l}{\textbf{Gurobi}} \\
$n$ & $d$ & opt & runtime [s] & gap & opt & runtime [s] & gap & opt & runtime [s] & gap\\\toprule
100 & 1.0 & 10/10 & 17, 104, 545 & - & 0/10 & - & 0.37 & 0/10 & - & 0.11\\
120 & 0.3 & 9/10 & 178, 539, 3189 & 0.18 & 10/10 & 52, 170, 930 & - & 10/10 & 10, 115, 821 & -\\
120 & 0.8 & 2/10 & 864, 1414, 1964 & 0.13 & 0/10 & - & 0.41 & 0/10 & - & 0.12\\
150 & 0.3 & 0/10 & - & 0.59 & 0/10 & - & 0.11 & 0/10 & - & 0.13\\
150 & 0.8 & 0/10 & - & 0.71 & 0/10 & - & 0.75 & 0/10 & - & 0.13\\
200 & 0.3 & 0/10 & - & 0.68 & 0/10 & - & 0.24 & 0/10 & - & 0.13\\
200 & 0.8 & 0/10 & - & 0.77 & 0/10 & - & 0.91 & 0/10 & - & 0.09\\
250 & 0.1 & 0/10 & - & 0.56 & 5/10 & 318, 543, 2913 & 0.04 & 6/10 & 941, 2452, 3541 & 0.03\\\bottomrule
\end{tabular}
\caption{Summary of results for \texttt{be*} instances from~\cite{billionnet2007using}. For each class of instances with $n$ variables and density $d$, we report (min,median,max) running times for instances that did not result in a timeout and the average gap for instances that resulted in a timeout. `-' means unavailable; either no instance was solved (runtime unavailable) or all instanced were solved (gap unavailable).}
\label{table:resultssummary_be}
\end{table}

\begin{table}[H]\footnotesize
\centering
\begin{tabular}{cccccccccc}
\multicolumn{1}{l}{\textbf{\texttt{bqp*} size}} & \multicolumn{3}{l}{\textbf{Solver}} & \multicolumn{3}{l}{\textbf{Xpress}} & \multicolumn{3}{l}{\textbf{Gurobi}} \\
$n$ & opt & runtime [s] & gap & opt & runtime [s] & gap & opt & runtime [s] & gap\\\toprule
50 & 10/10 & 0.1, 0.1, 1.1 & - & 10/10 & 0.0, 0.0, 0.0 & - & 10/10 & 0.0, 0.0, 0.0 & -\\
100 & 10/10 & 4.0, 11.3, 62.0 & - & 10/10 & 0.0, 0.1, 1.1 & - & 10/10 & 0.2, 0.2, 0.4 & -\\
250 & 0/10 & - & 0.57 & 5/10 & 372.9, 764.0, 1531.5 & 0.07 & 5/10 & 371.7, 2739.5, 3554.8 & 0.04\\
500 & 0/10 & - & 0.68 & 0/10 & - & 0.34 & 0/10 & - & 0.13\\\bottomrule
\end{tabular}
\caption{Summary of results for \texttt{bqp*} instances from~\cite{beasley1998heuristic}. For each class of instances with $n$ variables and density $d=0.1$, we report (min,median,max) running times for instances that did not result in a timeout and the average gap for instances that resulted in a timeout. `-' means unavailable; either no instance was solved (runtime unavailable) or all instanced were solved (gap unavailable).}
\label{table:resultssummary_bqp}
\end{table}

\begin{table}[H]\footnotesize
\centering
\begin{tabular}{cccccccccc}
\multicolumn{1}{l}{\textbf{\texttt{gka*}}} & \multicolumn{3}{l}{\textbf{Solver}} & \multicolumn{3}{l}{\textbf{Xpress}} & \multicolumn{3}{l}{\textbf{Gurobi}} \\
 & opt & runtime [s] & gap & opt & runtime [s] & gap & opt & runtime [s] & gap\\\toprule
& 28/45 & 0, 3, 562 & 0.84 & 31/45 & 0, 0, 2119 & 0.47 & 30/45 & 0, 0, 77 & 0.12\\\bottomrule
\end{tabular}
\caption{Summary of results for \texttt{gka*} instances from~\cite{glover1998adaptive}. We report (min,median,max) running times for instances that did not result in a timeout and the average gap for instances that resulted in a timeout.}
\label{table:resultssummary_gka}
\end{table}

\begin{table}[H]\footnotesize
\centering
\begin{tabular}{cccccccccc}
\multicolumn{1}{l}{\textbf{Ising problem}} & \multicolumn{3}{l}{\textbf{Solver}} & \multicolumn{3}{l}{\textbf{Xpress}} & \multicolumn{3}{l}{\textbf{Gurobi}} \\
Type & opt & runtime [s] & gap & opt & runtime [s] & gap & opt & runtime [s] & gap\\\toprule
dense, $\sigma=2.5$ & 15/15 & 0, 1, 33 & - & 6/15 & 42, 372, 3412 & 0.11 & 15/15 & 0, 0, 0 & -\\
dense, $\sigma=3.0$ & 15/15 & 0, 1, 5 & - & 11/15 & 19, 310, 2769 & 0.04 & 15/15 & 0, 0, 0 & -\\
2d & 5/9 & 0, 8, 1054 & 0.07 & 9/9 & 0, 1, 11 & - & 9/9 & 0, 0, 2 & -\\
3d & 2/9 & 300, 502, 704 & 0.21 & 9/9 & 1, 13, 550 & - & 9/9 & 0, 3, 20 & -\\\bottomrule
\end{tabular}
\caption{Summary of results for the Ising instances from~\cite{liers2004contributions}. For each class of instances, we report (min,median,max) running times for instances that did not result in a timeout and the average gap for instances that resulted in a timeout. `-' means unavailable; either no instance was solved (runtime unavailable) or all instanced were solved (gap unavailable).}
\label{table:resultssummary_ising}
\end{table}

\begin{table}[H]\footnotesize
\centering
\begin{tabular}{ccccccccccc}
\multicolumn{2}{l}{\textbf{Rudy problem}} & \multicolumn{3}{l}{\textbf{Solver}} & \multicolumn{3}{l}{\textbf{Xpress}} & \multicolumn{3}{l}{\textbf{Gurobi}} \\
$n$ & Type & opt & runtime [s] & gap & opt & runtime [s] & gap & opt & runtime [s] & gap\\\toprule
60 & g05 & 10/10 & 22, 50, 194 & - & 4/10 & 1700, 3135, 3471 & 0.03 & 10/10 & 50, 150, 203 & -\\
80 & g05 & 0/10 & - & 0.13 & 0/10 & - & 0.11 & 6/10 & 721, 2575, 3270 & 0.01\\
100 & g05 & 0/10 & - & 0.29 & 0/10 & - & 0.15 & 0/10 & - & 0.03\\\midrule
80 & pm1d & 10/10 & 8, 14, 40 & - & 0/10 & - & 0.56 & 0/10 & - & 0.20\\
100 & pm1d & 8/10 & 293, 904, 3332 & 0.18 & 0/10 & - & 0.69 & 0/10 & - & 0.20\\\midrule
80 & pm1s & 10/10 & 4, 6, 25 & - & 10/10 & 1, 3, 8 & - & 10/10 & 1, 2, 4 & -\\
100 & pm1s & 10/10 & 28, 234, 703 & - & 10/10 & 12, 22, 72 & - & 10/10 & 8, 19, 47 & -\\\midrule
100 & pw01 & 10/10 & 57, 127, 1487 & - & 10/10 & 10, 13, 28 & - & 10/10 & 6, 9, 12 & -\\
100 & pw05 & 0/10 & - & 0.26 & 0/10 & - & 0.12 & 0/10 & - & 0.04\\
100 & pw09 & 0/10 & - & 0.36 & 0/10 & - & 0.39 & 0/10 & - & 0.02\\\midrule
100 & w01 & 10/10 & 23, 56, 571 & - & 10/10 & 6, 14, 38 & - & 10/10 & 5, 6, 12 & -\\
100 & w05 & 10/10 & 103, 684, 2373 & - & 0/10 & - & 0.28 & 0/10 & - & 0.20\\
100 & w09 & 7/10 & 94, 482, 3112 & 0.18 & 0/10 & - & 0.59 & 0/10 & - & 0.20\\\bottomrule
\end{tabular}
\caption{Summary of results for rudy instances from~\cite{wiegele2007biq}. For each class of instances, we report (min,median,max) running times for instances that did not result in a timeout and the average gap for instances that resulted in a timeout. `-' means unavailable; either no instance was solved (runtime unavailable) or all instanced were solved (gap unavailable).}
\label{table:resultssummary_rudy}
\end{table}
\newpage
\subsection{Binary quadratic programming}

\begin{table}[H]
\centering
\begin{tabular}{lccc}
\textbf{Benchmark} & \textbf{Solver} & \textbf{Xpress} & \textbf{Gurobi}\\\toprule
be100.1.sparse & \textbf{-19412 (0.000)} & -18514 (0.347) & -19412 (0.092)\\
be100.2.sparse & \textbf{-17290 (0.000)} & -17132 (0.333) & -17290 (0.087)\\
be100.3.sparse & \textbf{-17565 (0.000)} & -16983 (0.314) & -17565 (0.091)\\
be100.4.sparse & \textbf{-19125 (0.000)} & -18538 (0.364) & -19125 (0.081)\\
be100.5.sparse & \textbf{-15868 (0.000)} & -15206 (0.339) & -15868 (0.106)\\
be100.6.sparse & \textbf{-17368 (0.000)} & -16462 (0.300) & -17368 (0.108)\\
be100.7.sparse & \textbf{-18629 (0.000)} & -17497 (0.429) & -18629 (0.104)\\
be100.8.sparse & \textbf{-18649 (0.000)} & -18199 (0.419) & -18649 (0.121)\\
be100.9.sparse & \textbf{-13294 (0.000)} & -11516 (0.436) & -13294 (0.142)\\
be100.10.sparse & \textbf{-15352 (0.000)} & -14543 (0.421) & -15352 (0.130)\\
\midrule
be120.3.1.sparse & \textbf{-13067 (0.000)} & \textbf{-13067 (0.000)} & \textbf{-13067 (0.000)}\\
be120.3.2.sparse & \textbf{-13046 (0.000)} & \textbf{-13046 (0.000)} & \textbf{-13046 (0.000)}\\
be120.3.3.sparse & \textbf{-12418 (0.000)} & \textbf{-12418 (0.000)} & \textbf{-12418 (0.000)}\\
be120.3.4.sparse & \textbf{-13867 (0.000)} & \textbf{-13867 (0.000)} & \textbf{-13867 (0.000)}\\
be120.3.5.sparse & \textbf{-11403 (0.000)} & \textbf{-11403 (0.000)} & \textbf{-11403 (0.000)}\\
be120.3.6.sparse & \textbf{-12915 (0.000)} & \textbf{-12915 (0.000)} & \textbf{-12915 (0.000)}\\
be120.3.7.sparse & \textbf{-14068 (0.000)} & \textbf{-14068 (0.000)} & \textbf{-14068 (0.000)}\\
be120.3.8.sparse & \textbf{-14701 (0.000)} & \textbf{-14701 (0.000)} & \textbf{-14701 (0.000)}\\
be120.3.9.sparse & -10458 (0.182) & \textbf{-10458 (0.000)} & \textbf{-10458 (0.000)}\\
be120.3.10.sparse & \textbf{-12201 (0.000)} & \textbf{-12201 (0.000)} & \textbf{-12201 (0.000)}\\
\midrule
be120.8.1.sparse & -18691 (0.140) & -8798 (0.715) & -18691 (0.133)\\
be120.8.2.sparse & -18827 (0.165) & -17484 (0.410) & -18827 (0.121)\\
be120.8.3.sparse & -19302 (0.123) & -18615 (0.344) & -19302 (0.117)\\
be120.8.4.sparse & \textbf{-20765 (0.000)} & -20288 (0.336) & -20765 (0.105)\\
be120.8.5.sparse & -20417 (0.091) & -18695 (0.347) & -20417 (0.094)\\
be120.8.6.sparse & -18482 (0.130) & -14570 (0.456) & -18482 (0.129)\\
be120.8.7.sparse & -22194 (0.117) & -21413 (0.396) & -22194 (0.112)\\
be120.8.8.sparse & -19534 (0.157) & -18268 (0.422) & -19534 (0.141)\\
be120.8.9.sparse & -18195 (0.146) & -17858 (0.317) & -18195 (0.136)\\
be120.8.10.sparse & \textbf{-19049 (0.000)} & -18292 (0.330) & -19049 (0.110)\\
\midrule
be150.3.1.sparse & -18889 (0.578) & -18782 (0.081) & -18889 (0.111)\\
be150.3.2.sparse & -17816 (0.600) & -17705 (0.103) & -17816 (0.128)\\
be150.3.3.sparse & -17314 (0.591) & -17314 (0.049) & -17314 (0.121)\\
be150.3.4.sparse & -19884 (0.550) & -19798 (0.063) & -19884 (0.100)\\
be150.3.5.sparse & -16817 (0.599) & -16810 (0.109) & -16817 (0.133)\\
be150.3.6.sparse & -16780 (0.605) & -16669 (0.129) & -16780 (0.141)\\
be150.3.7.sparse & -18001 (0.566) & -17836 (0.111) & -18001 (0.123)\\
be150.3.8.sparse & -18303 (0.589) & -18105 (0.145) & -18303 (0.128)\\
be150.3.9.sparse & -12838 (0.665) & -12546 (0.208) & -12838 (0.178)\\
be150.3.10.sparse & -17963 (0.582) & -17876 (0.114) & -17963 (0.125)\\
\midrule
be150.8.1.sparse & -27089 (0.717) & -11355 (0.771) & -27089 (0.122)\\
be150.8.2.sparse & -26779 (0.728) & -11024 (0.795) & -26779 (0.128)\\
be150.8.3.sparse & -29438 (0.698) & -17145 (0.672) & -29438 (0.123)\\
be150.8.4.sparse & -26911 (0.729) & -18917 (0.613) & -26911 (0.128)\\
be150.8.5.sparse & -28017 (0.708) & -12450 (0.744) & -28017 (0.123)\\
be150.8.6.sparse & -29221 (0.705) & -9960 (0.815) & -29221 (0.127)\\
be150.8.7.sparse & -31209 (0.672) & -27068 (0.554) & -31209 (0.118)\\
be150.8.8.sparse & -29730 (0.680) & -12907 (0.769) & -29730 (0.125)\\
be150.8.9.sparse & -25388 (0.740) & -1521 (0.967) & -25388 (0.146)\\
be150.8.10.sparse & -28374 (0.689) & -12091 (0.764) & -28374 (0.123)\\
\end{tabular}
\caption{Detailed results (final objective value and gap) for each \texttt{be*} instance and solver.}
\end{table}
\begin{table}[H]
\centering
\begin{tabular}{lccc}
\textbf{Benchmark} & \textbf{Solver} & \textbf{Xpress} & \textbf{Gurobi}\\\toprule
be200.3.1.sparse & -25453 (0.684) & -24550 (0.248) & -25453 (0.154)\\
be200.3.2.sparse & -25027 (0.687) & -24567 (0.220) & -25027 (0.131)\\
be200.3.3.sparse & -28023 (0.674) & -27227 (0.177) & -28023 (0.123)\\
be200.3.4.sparse & -27434 (0.666) & -26878 (0.186) & -27434 (0.132)\\
be200.3.5.sparse & -26355 (0.681) & -22778 (0.288) & -26355 (0.137)\\
be200.3.6.sparse & -26146 (0.677) & -21734 (0.336) & -26146 (0.133)\\
be200.3.7.sparse & -30483 (0.650) & -27491 (0.219) & -30483 (0.096)\\
be200.3.8.sparse & -27355 (0.666) & -26416 (0.212) & -27355 (0.137)\\
be200.3.9.sparse & -24683 (0.702) & -23444 (0.236) & -24683 (0.141)\\
be200.3.10.sparse & -23842 (0.696) & -22914 (0.251) & -23842 (0.153)\\\midrule
be200.8.1.sparse & -48534 (0.757) & -15347 (0.879) & -48534 (0.076)\\
be200.8.2.sparse & -40821 (0.787) & -9157 (0.912) & -40821 (0.120)\\
be200.8.3.sparse & -43207 (0.776) & -10323 (0.910) & -43207 (0.101)\\
be200.8.4.sparse & -43757 (0.778) & -7127 (0.934) & -43757 (0.084)\\
be200.8.5.sparse & -41482 (0.790) & -4866 (0.955) & -41482 (0.085)\\
be200.8.6.sparse & -49492 (0.756) & -11371 (0.902) & -49492 (0.053)\\
be200.8.7.sparse & -46828 (0.761) & -11168 (0.908) & -46828 (0.079)\\
be200.8.8.sparse & -44502 (0.777) & -16317 (0.870) & -44502 (0.092)\\
be200.8.9.sparse & -43241 (0.780) & -9728 (0.917) & -43241 (0.076)\\
be200.8.10.sparse & -42832 (0.783) & -8194 (0.925) & -42832 (0.095)\\
\midrule
be250.1.sparse & -24076 (0.540) & \textbf{-24076 (0.000)} & \textbf{-24076 (0.000)}\\
be250.2.sparse & -22540 (0.571) & -22540 (0.030) & \textbf{-22540 (0.000)}\\
be250.3.sparse & -22923 (0.555) & \textbf{-22923 (0.000)} & \textbf{-22923 (0.000)}\\
be250.4.sparse & -24649 (0.534) & \textbf{-24649 (0.000)} & \textbf{-24649 (0.000)}\\
be250.5.sparse & -21057 (0.590) & -21020 (0.062) & -21057 (0.048)\\
be250.6.sparse & -22735 (0.569) & -22709 (0.048) & -22735 (0.022)\\
be250.7.sparse & -24095 (0.557) & \textbf{-24095 (0.000)} & \textbf{-24095 (0.000)}\\
be250.8.sparse & -23801 (0.551) & -23801 (0.028) & \textbf{-23801 (0.000)}\\
be250.9.sparse & -20051 (0.593) & -20010 (0.038) & -20051 (0.017)\\
be250.10.sparse & -23159 (0.553) & \textbf{-23159 (0.000)} & -23159 (0.018)\\\bottomrule
\end{tabular}
\caption{(Continued:) Detailed results (final objective value and gap) for each \texttt{be*} instance and solver.}
\end{table}

\begin{table}[H]
\centering
\begin{tabular}{lccc}
\textbf{Benchmark} & \textbf{Solver} & \textbf{Xpress} & \textbf{Gurobi}\\\toprule
bqp50-1.sparse & \textbf{-2098 (0.000)} & \textbf{-2098 (0.000)} & \textbf{-2098 (0.000)}\\
bqp50-2.sparse & \textbf{-3702 (0.000)} & \textbf{-3702 (0.000)} & \textbf{-3702 (0.000)}\\
bqp50-3.sparse & \textbf{-4626 (0.000)} & \textbf{-4626 (0.000)} & \textbf{-4626 (0.000)}\\
bqp50-4.sparse & \textbf{-3544 (0.000)} & \textbf{-3544 (0.000)} & \textbf{-3544 (0.000)}\\
bqp50-5.sparse & \textbf{-4012 (0.000)} & \textbf{-4012 (0.000)} & \textbf{-4012 (0.000)}\\
bqp50-6.sparse & \textbf{-3693 (0.000)} & \textbf{-3693 (0.000)} & \textbf{-3693 (0.000)}\\
bqp50-7.sparse & \textbf{-4520 (0.000)} & \textbf{-4520 (0.000)} & \textbf{-4520 (0.000)}\\
bqp50-8.sparse & \textbf{-4216 (0.000)} & \textbf{-4216 (0.000)} & \textbf{-4216 (0.000)}\\
bqp50-9.sparse & \textbf{-3780 (0.000)} & \textbf{-3780 (0.000)} & \textbf{-3780 (0.000)}\\
bqp50-10.sparse & \textbf{-3507 (0.000)} & \textbf{-3507 (0.000)} & \textbf{-3507 (0.000)}\\
\midrule
bqp100-1.sparse & \textbf{-7970 (0.000)} & \textbf{-7970 (0.000)} & \textbf{-7970 (0.000)}\\
bqp100-2.sparse & \textbf{-11036 (0.000)} & \textbf{-11036 (0.000)} & \textbf{-11036 (0.000)}\\
bqp100-3.sparse & \textbf{-12723 (0.000)} & \textbf{-12723 (0.000)} & \textbf{-12723 (0.000)}\\
bqp100-4.sparse & \textbf{-10368 (0.000)} & \textbf{-10368 (0.000)} & \textbf{-10368 (0.000)}\\
bqp100-5.sparse & \textbf{-9083 (0.000)} & \textbf{-9083 (0.000)} & \textbf{-9083 (0.000)}\\
bqp100-6.sparse & \textbf{-10210 (0.000)} & \textbf{-10210 (0.000)} & \textbf{-10210 (0.000)}\\
bqp100-7.sparse & \textbf{-10125 (0.000)} & \textbf{-10125 (0.000)} & \textbf{-10125 (0.000)}\\
bqp100-8.sparse & \textbf{-11435 (0.000)} & \textbf{-11435 (0.000)} & \textbf{-11435 (0.000)}\\
bqp100-9.sparse & \textbf{-11455 (0.000)} & \textbf{-11455 (0.000)} & \textbf{-11455 (0.000)}\\
bqp100-10.sparse & \textbf{-12565 (0.000)} & \textbf{-12565 (0.000)} & \textbf{-12565 (0.000)}\\
\midrule
bqp250-1.sparse & -45607 (0.561) & \textbf{-45607 (0.000)} & \textbf{-45607 (0.000)}\\
bqp250-2.sparse & -44810 (0.557) & -44673 (0.045) & -44810 (0.016)\\
bqp250-3.sparse & -49037 (0.523) & \textbf{-49037 (0.000)} & \textbf{-49037 (0.000)}\\
bqp250-4.sparse & -41274 (0.582) & -41176 (0.035) & -41274 (0.023)\\
bqp250-5.sparse & -47961 (0.541) & \textbf{-47961 (0.000)} & \textbf{-47961 (0.000)}\\
bqp250-6.sparse & -41014 (0.594) & -40673 (0.079) & -41014 (0.049)\\
bqp250-7.sparse & -46757 (0.552) & \textbf{-46757 (0.000)} & \textbf{-46757 (0.000)}\\
bqp250-8.sparse & -35726 (0.626) & -35294 (0.112) & -35726 (0.073)\\
bqp250-9.sparse & -48916 (0.535) & \textbf{-48916 (0.000)} & \textbf{-48916 (0.000)}\\
bqp250-10.sparse & -40442 (0.582) & -40242 (0.062) & -40442 (0.036)\\
\midrule
bqp500-1.sparse & -116586 (0.696) & -100276 (0.362) & -116586 (0.137)\\
bqp500-2.sparse & -128339 (0.668) & -118452 (0.265) & -128339 (0.115)\\
bqp500-3.sparse & -130812 (0.669) & -113402 (0.350) & -130812 (0.115)\\
bqp500-4.sparse & -130097 (0.664) & -113947 (0.309) & -130097 (0.124)\\
bqp500-5.sparse & -125487 (0.678) & -108742 (0.344) & -125487 (0.131)\\
bqp500-6.sparse & -121772 (0.685) & -107398 (0.361) & -121772 (0.131)\\
bqp500-7.sparse & -122200 (0.689) & -105182 (0.376) & -122201 (0.136)\\
bqp500-8.sparse & -123559 (0.677) & -110131 (0.340) & -123559 (0.140)\\
bqp500-9.sparse & -120798 (0.687) & -102151 (0.365) & -120798 (0.131)\\
bqp500-10.sparse & -130619 (0.667) & -113383 (0.311) & -130619 (0.114)\\\bottomrule
\end{tabular}
\caption{Detailed results (final objective value and gap) for each \texttt{bqp*} instance and solver.}
\end{table}

\begin{table}[H]
\centering
\begin{tabular}{lccc}
\textbf{Benchmark} & \textbf{Solver} & \textbf{Xpress} & \textbf{Gurobi}\\\toprule
gka1a.sparse & \textbf{-3414 (0.000)} & \textbf{-3414 (0.000)} & \textbf{-3414 (0.000)}\\
gka1b.sparse & \textbf{-133 (0.000)} & \textbf{-133 (0.000)} & \textbf{-133 (0.000)}\\
gka1c.sparse & \textbf{-5058 (0.000)} & \textbf{-5058 (0.000)} & \textbf{-5058 (0.000)}\\
gka1d.sparse & \textbf{-6333 (0.000)} & \textbf{-6333 (0.000)} & \textbf{-6333 (0.000)}\\
gka1e.sparse & -16464 (0.485) & \textbf{-16464 (0.000)} & \textbf{-16464 (0.000)}\\
gka1f.sparse & -61194 (0.693) & -51659 (0.348) & -61194 (0.133)\\
\midrule
gka2a.sparse & \textbf{-6063 (0.000)} & \textbf{-6063 (0.000)} & \textbf{-6063 (0.000)}\\
gka2b.sparse & \textbf{-121 (0.000)} & \textbf{-121 (0.000)} & \textbf{-121 (0.000)}\\
gka2c.sparse & \textbf{-6213 (0.000)} & \textbf{-6213 (0.000)} & \textbf{-6213 (0.000)}\\
gka2d.sparse & \textbf{-6579 (0.000)} & \textbf{-6579 (0.000)} & \textbf{-6579 (0.000)}\\
gka2e.sparse & -23395 (0.608) & -23199 (0.104) & -23395 (0.068)\\
gka2f.sparse & -100161 (0.780) & -23590 (0.908) & -100161 (0.131)\\
\midrule
gka3a.sparse & \textbf{-6037 (0.000)} & \textbf{-6037 (0.000)} & \textbf{-6037 (0.000)}\\
gka3b.sparse & \textbf{-118 (0.000)} & \textbf{-118 (0.000)} & \textbf{-118 (0.000)}\\
gka3c.sparse & \textbf{-6665 (0.000)} & \textbf{-6665 (0.000)} & \textbf{-6665 (0.000)}\\
gka3d.sparse & \textbf{-9261 (0.000)} & \textbf{-9261 (0.000)} & \textbf{-9261 (0.000)}\\
gka3e.sparse & -25243 (0.688) & -24714 (0.199) & -25243 (0.137)\\
gka3f.sparse & -138035 (0.840) & -37112 (0.943) & -138035 (0.140)\\
\midrule
gka4a.sparse & \textbf{-8598 (0.000)} & \textbf{-8598 (0.000)} & \textbf{-8598 (0.000)}\\
gka4b.sparse & -129 (0.991) & \textbf{-129 (0.000)} & \textbf{-129 (0.000)}\\
gka4c.sparse & \textbf{-7398 (0.000)} & \textbf{-7398 (0.000)} & \textbf{-7398 (0.000)}\\
gka4d.sparse & \textbf{-10727 (0.000)} & \textbf{-10727 (0.000)} & \textbf{-10727 (0.000)}\\
gka4e.sparse & -35594 (0.668) & -27862 (0.398) & -35594 (0.099)\\
gka4f.sparse & -172771 (0.864) & -52836 (0.953) & -172771 (0.135)\\
\midrule
gka5a.sparse & \textbf{-5737 (0.000)} & \textbf{-5737 (0.000)} & \textbf{-5737 (0.000)}\\
gka5b.sparse & -150 (0.995) & \textbf{-150 (0.000)} & \textbf{-150 (0.000)}\\
gka5c.sparse & \textbf{-7362 (0.000)} & \textbf{-7362 (0.000)} & \textbf{-7362 (0.000)}\\
gka5d.sparse & \textbf{-11626 (0.000)} & -11589 (0.114) & -11626 (0.124)\\
gka5e.sparse & -35154 (0.727) & -29961 (0.526) & -35154 (0.105)\\
gka5f.sparse & -190507 (0.885) & -47758 (0.968) & -190507 (0.137)\\
\midrule
gka6a.sparse & \textbf{-3980 (0.000)} & \textbf{-3980 (0.000)} & \textbf{-3980 (0.000)}\\
gka6b.sparse & -146 (0.997) & \textbf{-146 (0.000)} & \textbf{-146 (0.000)}\\
gka6c.sparse & \textbf{-5824 (0.000)} & \textbf{-5824 (0.000)} & \textbf{-5824 (0.000)}\\
gka6d.sparse & \textbf{-14207 (0.000)} & \textbf{-14207 (0.000)} & -14207 (0.099)\\
\midrule
gka7a.sparse & \textbf{-4541 (0.000)} & \textbf{-4541 (0.000)} & \textbf{-4541 (0.000)}\\
gka7b.sparse & -160 (0.998) & \textbf{-160 (0.000)} & \textbf{-160 (0.000)}\\
gka7c.sparse & \textbf{-7225 (0.000)} & \textbf{-7225 (0.000)} & \textbf{-7225 (0.000)}\\
gka7d.sparse & \textbf{-14476 (0.000)} & -14241 (0.230) & -14476 (0.114)\\
\midrule
gka8a.sparse & \textbf{-11109 (0.000)} & \textbf{-11109 (0.000)} & \textbf{-11109 (0.000)}\\
gka8b.sparse & -145 (0.999) & \textbf{-145 (0.000)} & \textbf{-145 (0.000)}\\
gka8d.sparse & \textbf{-16352 (0.000)} & -16045 (0.229) & -16352 (0.088)\\
\midrule
gka9b.sparse & -137 (0.999) & \textbf{-137 (0.000)} & \textbf{-137 (0.000)}\\
gka9d.sparse & \textbf{-15656 (0.000)} & -15199 (0.308) & -15656 (0.123)\\
\midrule
gka10b.sparse & -154 (0.999) & \textbf{-154 (0.000)} & \textbf{-154 (0.000)}\\
gka10d.sparse & \textbf{-19102 (0.000)} & -18399 (0.413) & -19102 (0.091)\\\bottomrule
\end{tabular}
\caption{Detailed results (final objective value and gap) for each \texttt{gka*} instance and solver.}
\end{table}

\newpage
\subsection{MaxCut}

\begin{table}[H]
\centering
\begin{tabular}{lccc}
\textbf{Benchmark} & \textbf{Solver} & \textbf{Xpress} & \textbf{Gurobi}\\\toprule
ising2.5-100\_5555 & \textbf{2460049 (0.000)} & \textbf{2460049 (0.000)} & \textbf{2460049 (0.000)}\\
ising2.5-100\_6666 & \textbf{2031217 (0.000)} & \textbf{2031217 (0.000)} & \textbf{2031217 (0.000)}\\
ising2.5-100\_7777 & \textbf{3363230 (0.000)} & \textbf{3363230 (0.000)} & \textbf{3363230 (0.000)}\\
ising2.5-150\_5555 & \textbf{4363532 (0.000)} & \textbf{4363532 (0.000)} & \textbf{4363532 (0.000)}\\
ising2.5-150\_6666 & \textbf{4057153 (0.000)} & \textbf{4057153 (0.000)} & \textbf{4057153 (0.000)}\\
ising2.5-150\_7777 & \textbf{4243269 (0.000)} & \textbf{4243269 (0.000)} & \textbf{4243269 (0.000)}\\
ising2.5-200\_5555 & \textbf{6294701 (0.000)} & 6278964 (0.072) & \textbf{6294701 (0.000)}\\
ising2.5-200\_6666 & \textbf{6795365 (0.000)} & 6778180 (0.090) & \textbf{6795365 (0.000)}\\
ising2.5-200\_7777 & \textbf{5568272 (0.000)} & 5561673 (0.081) & \textbf{5568272 (0.000)}\\
ising2.5-250\_5555 & \textbf{7919449 (0.000)} & 7845137 (0.087) & \textbf{7919449 (0.000)}\\
ising2.5-250\_6666 & \textbf{6925717 (0.000)} & 6856949 (0.124) & \textbf{6925717 (0.000)}\\
ising2.5-250\_7777 & \textbf{6596797 (0.000)} & 6550580 (0.100) & \textbf{6596797 (0.000)}\\
ising2.5-300\_5555 & \textbf{8579363 (0.000)} & 8412040 (0.157) & \textbf{8579363 (0.000)}\\
ising2.5-300\_6666 & \textbf{9102033 (0.000)} & 8973099 (0.106) & \textbf{9102033 (0.000)}\\
ising2.5-300\_7777 & \textbf{8323804 (0.000)} & 8199246 (0.153) & \textbf{8323804 (0.000)}\\
\midrule
ising3.0-100\_5555 & \textbf{2448189 (0.000)} & \textbf{2448189 (0.000)} & \textbf{2448189 (0.000)}\\
ising3.0-100\_6666 & \textbf{1984099 (0.000)} & \textbf{1984099 (0.000)} & \textbf{1984099 (0.000)}\\
ising3.0-100\_7777 & \textbf{3335814 (0.000)} & \textbf{3335814 (0.000)} & \textbf{3335814 (0.000)}\\
ising3.0-150\_5555 & \textbf{4279261 (0.000)} & \textbf{4279261 (0.000)} & \textbf{4279261 (0.000)}\\
ising3.0-150\_6666 & \textbf{3949317 (0.000)} & \textbf{3949317 (0.000)} & \textbf{3949317 (0.000)}\\
ising3.0-150\_7777 & \textbf{4211158 (0.000)} & \textbf{4211158 (0.000)} & \textbf{4211158 (0.000)}\\
ising3.0-200\_5555 & \textbf{6215531 (0.000)} & \textbf{6215531 (0.000)} & \textbf{6215531 (0.000)}\\
ising3.0-200\_6666 & \textbf{6756263 (0.000)} & \textbf{6756263 (0.000)} & \textbf{6756263 (0.000)}\\
ising3.0-200\_7777 & \textbf{5560824 (0.000)} & \textbf{5560824 (0.000)} & \textbf{5560824 (0.000)}\\
ising3.0-250\_5555 & \textbf{7823791 (0.000)} & 7823791 (0.039) & \textbf{7823791 (0.000)}\\
ising3.0-250\_6666 & \textbf{6903351 (0.000)} & \textbf{6903351 (0.000)} & \textbf{6903351 (0.000)}\\
ising3.0-250\_7777 & \textbf{6418276 (0.000)} & \textbf{6418276 (0.000)} & \textbf{6418276 (0.000)}\\
ising3.0-300\_5555 & \textbf{8493173 (0.000)} & 8489933 (0.039) & \textbf{8493173 (0.000)}\\
ising3.0-300\_6666 & \textbf{8915110 (0.000)} & 8914176 (0.040) & \textbf{8915110 (0.000)}\\
ising3.0-300\_7777 & \textbf{8242904 (0.000)} & 8241481 (0.040) & \textbf{8242904 (0.000)}\\
\midrule
t2g10\_5555 & \textbf{6049461 (0.000)} & \textbf{6049461 (0.000)} & \textbf{6049461 (0.000)}\\
t2g10\_6666 & \textbf{5757868 (0.000)} & \textbf{5757868 (0.000)} & \textbf{5757868 (0.000)}\\
t2g10\_7777 & \textbf{6509837 (0.000)} & \textbf{6509837 (0.000)} & \textbf{6509837 (0.000)}\\
t2g15\_5555 & \textbf{15051133 (0.000)} & \textbf{15051133 (0.000)} & \textbf{15051133 (0.000)}\\
t2g15\_6666 & 15697581 (0.019) & \textbf{15763716 (0.000)} & \textbf{15763716 (0.000)}\\
t2g15\_7777 & \textbf{15269399 (0.000)} & \textbf{15269399 (0.000)} & \textbf{15269399 (0.000)}\\
t2g20\_5555 & 24750563 (0.110) & \textbf{24838942 (0.000)} & \textbf{24838942 (0.000)}\\
t2g20\_6666 & 29185429 (0.066) & \textbf{29290570 (0.000)} & \textbf{29290570 (0.000)}\\
t2g20\_7777 & 28193799 (0.088) & \textbf{28349398 (0.000)} & \textbf{28349398 (0.000)}\\
\midrule
t3g5\_5555 & \textbf{10933215 (0.000)} & \textbf{10933215 (0.000)} & \textbf{10933215 (0.000)}\\
t3g5\_6666 & 11582216 (0.062) & \textbf{11582216 (0.000)} & \textbf{11582216 (0.000)}\\
t3g5\_7777 & \textbf{11552046 (0.000)} & \textbf{11552046 (0.000)} & \textbf{11552046 (0.000)}\\
t3g6\_5555 & 16893302 (0.218) & \textbf{17434469 (0.000)} & \textbf{17434469 (0.000)}\\
t3g6\_6666 & 19820954 (0.185) & \textbf{20217380 (0.000)} & \textbf{20217380 (0.000)}\\
t3g6\_7777 & 19159401 (0.194) & \textbf{19475011 (0.000)} & \textbf{19475011 (0.000)}\\
t3g7\_5555 & 27350066 (0.285) & \textbf{28302918 (0.000)} & \textbf{28302918 (0.000)}\\
t3g7\_6666 & 32531022 (0.239) & \textbf{33611981 (0.000)} & \textbf{33611981 (0.000)}\\
t3g7\_7777 & 28485184 (0.278) & \textbf{29118445 (0.000)} & \textbf{29118445 (0.000)}\\\bottomrule
\end{tabular}
\caption{Detailed results (final objective value and gap) for each Ising instance and solver.}
\end{table}

\begin{table}[H]
\centering
\begin{tabular}{lccc}
\textbf{Benchmark} & \textbf{Solver} & \textbf{Xpress} & \textbf{Gurobi}\\\toprule
g05\_60.0 & \textbf{536 (0.000)} & \textbf{536 (0.000)} & \textbf{536 (0.000)}\\
g05\_60.1 & \textbf{532 (0.000)} & 532 (0.020) & \textbf{532 (0.000)}\\
g05\_60.2 & \textbf{529 (0.000)} & 529 (0.024) & \textbf{529 (0.000)}\\
g05\_60.3 & \textbf{538 (0.000)} & \textbf{538 (0.000)} & \textbf{538 (0.000)}\\
g05\_60.4 & \textbf{527 (0.000)} & 527 (0.035) & \textbf{527 (0.000)}\\
g05\_60.5 & \textbf{533 (0.000)} & \textbf{533 (0.000)} & \textbf{533 (0.000)}\\
g05\_60.6 & \textbf{531 (0.000)} & 531 (0.025) & \textbf{531 (0.000)}\\
g05\_60.7 & \textbf{535 (0.000)} & \textbf{535 (0.000)} & \textbf{535 (0.000)}\\
g05\_60.8 & \textbf{530 (0.000)} & 530 (0.034) & \textbf{530 (0.000)}\\
g05\_60.9 & \textbf{533 (0.000)} & 533 (0.032) & \textbf{533 (0.000)}\\
\midrule
g05\_80.0 & 929 (0.133) & 917 (0.107) & \textbf{929 (0.000)}\\
g05\_80.1 & 941 (0.124) & 935 (0.078) & \textbf{941 (0.000)}\\
g05\_80.2 & 934 (0.135) & 920 (0.108) & \textbf{934 (0.000)}\\
g05\_80.3 & 923 (0.145) & 915 (0.106) & 923 (0.011)\\
g05\_80.4 & 932 (0.131) & 924 (0.107) & \textbf{932 (0.000)}\\
g05\_80.5 & 926 (0.138) & 917 (0.096) & 926 (0.005)\\
g05\_80.6 & 929 (0.135) & 918 (0.106) & \textbf{929 (0.000)}\\
g05\_80.7 & 929 (0.127) & 919 (0.106) & \textbf{929 (0.000)}\\
g05\_80.8 & 925 (0.129) & 914 (0.119) & 925 (0.006)\\
g05\_80.9 & 923 (0.139) & 912 (0.118) & 923 (0.009)\\
\midrule
g05\_100.0 & 1430 (0.301) & 1411 (0.140) & 1430 (0.035)\\
g05\_100.1 & 1425 (0.291) & 1392 (0.150) & 1425 (0.037)\\
g05\_100.2 & 1432 (0.282) & 1383 (0.157) & 1432 (0.030)\\
g05\_100.3 & 1424 (0.298) & 1390 (0.151) & 1424 (0.031)\\
g05\_100.4 & 1440 (0.292) & 1414 (0.136) & 1440 (0.030)\\
g05\_100.5 & 1436 (0.303) & 1394 (0.149) & 1436 (0.028)\\
g05\_100.6 & 1434 (0.296) & 1403 (0.144) & 1434 (0.028)\\
g05\_100.7 & 1431 (0.294) & 1392 (0.147) & 1431 (0.035)\\
g05\_100.8 & 1432 (0.280) & 1415 (0.137) & 1432 (0.030)\\
g05\_100.9 & 1430 (0.297) & 1396 (0.148) & 1430 (0.031)\\
\midrule
pm1d\_80.0 & \textbf{227 (0.000)} & 199 (0.566) & 227 (0.225)\\
pm1d\_80.1 & \textbf{245 (0.000)} & 217 (0.554) & 245 (0.202)\\
pm1d\_80.2 & \textbf{284 (0.000)} & 252 (0.522) & 284 (0.184)\\
pm1d\_80.3 & \textbf{291 (0.000)} & 260 (0.502) & 291 (0.176)\\
pm1d\_80.4 & \textbf{251 (0.000)} & 210 (0.569) & 251 (0.196)\\
pm1d\_80.5 & \textbf{242 (0.000)} & 209 (0.573) & 242 (0.214)\\
pm1d\_80.6 & \textbf{205 (0.000)} & 184 (0.584) & 205 (0.229)\\
pm1d\_80.7 & \textbf{249 (0.000)} & 208 (0.577) & 249 (0.192)\\
pm1d\_80.8 & \textbf{293 (0.000)} & 248 (0.538) & 293 (0.182)\\
pm1d\_80.9 & \textbf{258 (0.000)} & 210 (0.575) & 258 (0.160)\\
\midrule
pm1d\_100.0 & 340 (0.194) & 212 (0.728) & 340 (0.218)\\
pm1d\_100.1 & \textbf{324 (0.000)} & 188 (0.761) & 324 (0.225)\\
pm1d\_100.2 & 389 (0.167) & 232 (0.728) & 389 (0.181)\\
pm1d\_100.3 & \textbf{400 (0.000)} & 284 (0.657) & 400 (0.187)\\
pm1d\_100.4 & \textbf{363 (0.000)} & 218 (0.726) & 363 (0.214)\\
pm1d\_100.5 & \textbf{441 (0.000)} & 253 (0.712) & 441 (0.188)\\
pm1d\_100.6 & \textbf{367 (0.000)} & 253 (0.676) & 367 (0.192)\\
pm1d\_100.7 & \textbf{361 (0.000)} & 216 (0.728) & 361 (0.220)\\
pm1d\_100.8 & \textbf{385 (0.000)} & 309 (0.612) & 385 (0.184)\\
pm1d\_100.9 & \textbf{405 (0.000)} & 337 (0.589) & 405 (0.182)\\
\bottomrule
\end{tabular}
\caption{Detailed results (final objective value and gap) for each rudy instance and solver.}
\end{table}

\begin{table}[H]
\centering
\begin{tabular}{lccc}
\textbf{Benchmark} & \textbf{Solver} & \textbf{Xpress} & \textbf{Gurobi}\\\toprule
pm1s\_80.0 & \textbf{79 (0.000)} & \textbf{79 (0.000)} & \textbf{79 (0.000)}\\
pm1s\_80.1 & \textbf{85 (0.000)} & \textbf{85 (0.000)} & \textbf{85 (0.000)}\\
pm1s\_80.2 & \textbf{82 (0.000)} & \textbf{82 (0.000)} & \textbf{82 (0.000)}\\
pm1s\_80.3 & \textbf{81 (0.000)} & \textbf{81 (0.000)} & \textbf{81 (0.000)}\\
pm1s\_80.4 & \textbf{70 (0.000)} & \textbf{70 (0.000)} & \textbf{70 (0.000)}\\
pm1s\_80.5 & \textbf{87 (0.000)} & \textbf{87 (0.000)} & \textbf{87 (0.000)}\\
pm1s\_80.6 & \textbf{73 (0.000)} & \textbf{73 (0.000)} & \textbf{73 (0.000)}\\
pm1s\_80.7 & \textbf{83 (0.000)} & \textbf{83 (0.000)} & \textbf{83 (0.000)}\\
pm1s\_80.8 & \textbf{81 (0.000)} & \textbf{81 (0.000)} & \textbf{81 (0.000)}\\
pm1s\_80.9 & \textbf{70 (0.000)} & \textbf{70 (0.000)} & \textbf{70 (0.000)}\\
\midrule
pm1s\_100.0 & \textbf{127 (0.000)} & \textbf{127 (0.000)} & \textbf{127 (0.000)}\\
pm1s\_100.1 & \textbf{126 (0.000)} & \textbf{126 (0.000)} & \textbf{126 (0.000)}\\
pm1s\_100.2 & \textbf{125 (0.000)} & \textbf{125 (0.000)} & \textbf{125 (0.000)}\\
pm1s\_100.3 & \textbf{111 (0.000)} & \textbf{111 (0.000)} & \textbf{111 (0.000)}\\
pm1s\_100.4 & \textbf{128 (0.000)} & \textbf{128 (0.000)} & \textbf{128 (0.000)}\\
pm1s\_100.5 & \textbf{128 (0.000)} & \textbf{128 (0.000)} & \textbf{128 (0.000)}\\
pm1s\_100.6 & \textbf{122 (0.000)} & \textbf{122 (0.000)} & \textbf{122 (0.000)}\\
pm1s\_100.7 & \textbf{112 (0.000)} & \textbf{112 (0.000)} & \textbf{112 (0.000)}\\
pm1s\_100.8 & \textbf{120 (0.000)} & \textbf{120 (0.000)} & \textbf{120 (0.000)}\\
pm1s\_100.9 & \textbf{127 (0.000)} & \textbf{127 (0.000)} & \textbf{127 (0.000)}\\
\midrule
pw01\_100.0 & \textbf{2019 (0.000)} & \textbf{2019 (0.000)} & \textbf{2019 (0.000)}\\
pw01\_100.1 & \textbf{2060 (0.000)} & \textbf{2060 (0.000)} & \textbf{2060 (0.000)}\\
pw01\_100.2 & \textbf{2032 (0.000)} & \textbf{2032 (0.000)} & \textbf{2032 (0.000)}\\
pw01\_100.3 & \textbf{2067 (0.000)} & \textbf{2067 (0.000)} & \textbf{2067 (0.000)}\\
pw01\_100.4 & \textbf{2039 (0.000)} & \textbf{2039 (0.000)} & \textbf{2039 (0.000)}\\
pw01\_100.5 & \textbf{2108 (0.000)} & \textbf{2108 (0.000)} & \textbf{2108 (0.000)}\\
pw01\_100.6 & \textbf{2032 (0.000)} & \textbf{2032 (0.000)} & \textbf{2032 (0.000)}\\
pw01\_100.7 & \textbf{2074 (0.000)} & \textbf{2074 (0.000)} & \textbf{2074 (0.000)}\\
pw01\_100.8 & \textbf{2022 (0.000)} & \textbf{2022 (0.000)} & \textbf{2022 (0.000)}\\
pw01\_100.9 & \textbf{2005 (0.000)} & \textbf{2005 (0.000)} & \textbf{2005 (0.000)}\\
\midrule
pw05\_100.0 & 8190 (0.255) & 8043 (0.120) & 8190 (0.042)\\
pw05\_100.1 & 8045 (0.270) & 7894 (0.119) & 8045 (0.037)\\
pw05\_100.2 & 8039 (0.245) & 7939 (0.112) & 8039 (0.040)\\
pw05\_100.3 & 8139 (0.274) & 7897 (0.129) & 8139 (0.034)\\
pw05\_100.4 & 8125 (0.273) & 8009 (0.117) & 8125 (0.042)\\
pw05\_100.5 & 8169 (0.233) & 7980 (0.122) & 8169 (0.038)\\
pw05\_100.6 & 8217 (0.258) & 8076 (0.118) & 8217 (0.043)\\
pw05\_100.7 & 8249 (0.234) & 8086 (0.116) & 8249 (0.040)\\
pw05\_100.8 & 8199 (0.259) & 8080 (0.111) & 8199 (0.032)\\
pw05\_100.9 & 8099 (0.270) & 7927 (0.123) & 8099 (0.042)\\
\midrule
pw09\_100.0 & 13585 (0.365) & 13065 (0.202) & 13585 (0.021)\\
pw09\_100.1 & 13417 (0.362) & 13001 (0.194) & 13417 (0.023)\\
pw09\_100.2 & 13461 (0.363) & 2440 (0.850) & 13461 (0.022)\\
pw09\_100.3 & 13656 (0.365) & 13171 (0.193) & 13656 (0.024)\\
pw09\_100.4 & 13514 (0.363) & 12992 (0.202) & 13514 (0.020)\\
pw09\_100.5 & 13574 (0.364) & 2422 (0.852) & 13574 (0.024)\\
pw09\_100.6 & 13640 (0.362) & 3381 (0.794) & 13640 (0.018)\\
pw09\_100.7 & 13501 (0.364) & 13157 (0.192) & 13501 (0.022)\\
pw09\_100.8 & 13593 (0.363) & 13124 (0.197) & 13593 (0.022)\\
pw09\_100.9 & 13658 (0.359) & 13360 (0.182) & 13658 (0.022)\\
\bottomrule
\end{tabular}
\caption{(Continued 2/3:) Detailed results (final objective value and gap) for each rudy instance and solver.}
\end{table}

\begin{table}[H]
\centering
\begin{tabular}{lccc}
\textbf{Benchmark} & \textbf{Solver} & \textbf{Xpress} & \textbf{Gurobi}\\\toprule
w01\_100.0 & \textbf{651 (0.000)} & \textbf{651 (0.000)} & \textbf{651 (0.000)}\\
w01\_100.1 & \textbf{719 (0.000)} & \textbf{719 (0.000)} & \textbf{719 (0.000)}\\
w01\_100.2 & \textbf{676 (0.000)} & \textbf{676 (0.000)} & \textbf{676 (0.000)}\\
w01\_100.3 & \textbf{813 (0.000)} & \textbf{813 (0.000)} & \textbf{813 (0.000)}\\
w01\_100.4 & \textbf{668 (0.000)} & \textbf{668 (0.000)} & \textbf{668 (0.000)}\\
w01\_100.5 & \textbf{643 (0.000)} & \textbf{643 (0.000)} & \textbf{643 (0.000)}\\
w01\_100.6 & \textbf{654 (0.000)} & \textbf{654 (0.000)} & \textbf{654 (0.000)}\\
w01\_100.7 & \textbf{725 (0.000)} & \textbf{725 (0.000)} & \textbf{725 (0.000)}\\
w01\_100.8 & \textbf{721 (0.000)} & \textbf{721 (0.000)} & \textbf{721 (0.000)}\\
w01\_100.9 & \textbf{729 (0.000)} & \textbf{729 (0.000)} & \textbf{729 (0.000)}\\
\midrule
w05\_100.0 & \textbf{1646 (0.000)} & 1596 (0.280) & 1646 (0.209)\\
w05\_100.1 & \textbf{1606 (0.000)} & 1543 (0.290) & 1606 (0.185)\\
w05\_100.2 & \textbf{1902 (0.000)} & 1869 (0.231) & 1902 (0.188)\\
w05\_100.3 & \textbf{1627 (0.000)} & 1515 (0.308) & 1627 (0.197)\\
w05\_100.4 & \textbf{1546 (0.000)} & 1437 (0.319) & 1546 (0.225)\\
w05\_100.5 & \textbf{1581 (0.000)} & 1524 (0.294) & 1581 (0.204)\\
w05\_100.6 & \textbf{1479 (0.000)} & 1397 (0.312) & 1479 (0.212)\\
w05\_100.7 & \textbf{1987 (0.000)} & 1920 (0.222) & 1987 (0.171)\\
w05\_100.8 & \textbf{1311 (0.000)} & 1263 (0.329) & 1311 (0.235)\\
w05\_100.9 & \textbf{1752 (0.000)} & 1727 (0.215) & 1752 (0.184)\\
\midrule
w09\_100.0 & \textbf{2121 (0.000)} & 1733 (0.543) & 2121 (0.195)\\
w09\_100.1 & 2096 (0.205) & 1296 (0.661) & 2096 (0.225)\\
w09\_100.2 & 2738 (0.103) & 1517 (0.660) & 2738 (0.177)\\
w09\_100.3 & 1990 (0.220) & 1490 (0.601) & 1990 (0.202)\\
w09\_100.4 & \textbf{2033 (0.000)} & 1747 (0.522) & 2033 (0.217)\\
w09\_100.5 & \textbf{2433 (0.000)} & 1863 (0.528) & 2433 (0.170)\\
w09\_100.6 & \textbf{2220 (0.000)} & 1995 (0.474) & 2220 (0.214)\\
w09\_100.7 & \textbf{2252 (0.000)} & 1966 (0.493) & 2252 (0.195)\\
w09\_100.8 & \textbf{1843 (0.000)} & 1035 (0.693) & 1843 (0.240)\\
w09\_100.9 & \textbf{2043 (0.000)} & 1142 (0.689) & 2043 (0.207)\\\bottomrule
\end{tabular}
\caption{(Continued 3/3:) Detailed results (final objective value and gap) for each rudy instance and solver.}
\end{table}

\end{document}